\documentclass[11pt,reqno]{amsart}

\usepackage[mathscr]{eucal}
\usepackage[all]{xy}
\usepackage{mathrsfs}
\usepackage{xypic}
\usepackage{amsfonts}
\usepackage{amsmath}
\usepackage{amsthm}
\usepackage{amssymb}
\usepackage{latexsym}
\usepackage{eufrak}

\theoremstyle{plain}
\newtheorem{definition}[equation]{Definition}
\newtheorem{corollary}[equation]{Corollary}
\newtheorem{lemma}[equation]{Lemma}
\newtheorem{proposition}[equation]{Proposition}
\newtheorem{theorem}[equation]{Theorem}

\theoremstyle{remark}
\newtheorem{remark}[equation]{Remark}
\newtheorem{example}[equation]{Example}
\newtheorem{notation}[equation]{Notation}

\numberwithin{equation}{subsection}

\makeatletter
\newcommand{\subsec}{\@startsection{subsection}{2}{0pt}{-3ex
plus -1ex minus -0.2ex}{-2mm plus -0pt minus
-2pt}{\normalfont\bfseries}} \makeatother

\newcommand{\erem}{\hphantom{.}\hfill$\lozenge$\end{remark}}

\newcommand{\Lmod}[1]{#1\text{-}{\mathsf{mod}}}

\newcommand{\bimod}[1]{#1\text{-}{\sf{bimod}}}

\newcommand{\hdot}{{\:\raisebox{2pt}{\text{\circle*{1.5}}}}}
\newcommand{\idot}{{\:\raisebox{2pt}{\text{\circle*{1.5}}}}}

\DeclareMathOperator{\Ext}{\mathrm{Ext}}

\DeclareMathOperator{\Tr}{\mathrm{Tr}}
\DeclareMathOperator{\gr}{\mathrm{gr}}
\DeclareMathOperator{\Rep}{\mathrm{Rep}}

\DeclareMathOperator{\im}{\mathrm{Im}}
\DeclareMathOperator{\Coker}{\mathrm{Coker}}
\DeclareMathOperator{\sym}{{\mathrm{Sym}}}
\DeclareMathOperator{\ssym}{{\mathrm{Sym}}}

\DeclareMathOperator{\La}{{\mathsf{\Lambda}}}
\DeclareMathOperator{\Id}{{\mathrm{Id}}}

\DeclareMathOperator{\lstl}{%
\hbox{\:\,\vbox to 0pt{\vss{\vbox to 0pt
{\vss\hbox to 0 pt {\Huge\hss$\circ$\hss}\vss}%
\vskip -14.2pt%
\vbox to 0pt {\vss\hbox to 0 pt {\hss$\ltimes$\hss}\vss}%
\vskip 3pt}}}\:\,}

\newcommand{\dis}{\displaystyle}

\newcommand{\beq}{\begin{equation}\label}
\newcommand{\eeq}{\end{equation}}

\newcommand{\iso}{{\;\stackrel{_\sim}{\to}\;}}
\newcommand{\cd}{\!\cdot\!}

\def\ccirc{{{}_{\,{}^{^\circ}}}}

\DeclareMathOperator{\End}{\mathrm{End}}

\DeclareMathOperator{\Hom}{\mathrm{Hom}}

\newcommand{\Ker}{\text{Ker\,}}

\newcommand{\z}{^\flat }

\newcommand{\Om}{\Omega }
\newcommand{\DR}{\overline{\operatorname{DR}}}

\newcommand{\HC}{{\overline{HC}}}

\newcommand{\RR}{{\mathsf R}}
\newcommand{\bc}{{\mathbf{c}}}
\newcommand{\bm}{{\mathbf{m}}}

\renewcommand{\O}{\mathcal{O}}

\newcommand{\G}{\Gamma }

\newcommand{\wh}{\widehat}

\newcommand{\ta}{{A\z}}

\newcommand{\sminus}{\smallsetminus}
\newcommand{\mto}{\mapsto}

\newcommand{\inv}{^{-1}}
\newcommand{\vi}{${\en\sf {(i)}}\;$}
\newcommand{\vii}{${\;\sf {(ii)}}\;$}
\newcommand{\viii}{${\sf {(iii)}}\;$}

\newcommand{\sset}{\subset}

\newcommand{\into}{{}^{\,}\hookrightarrow^{\,}}
\newcommand{\too}{\,\longrightarrow\,}
\newcommand{\onto}{\twoheadrightarrow}
\newcommand{\en}{\enspace }
\renewcommand{\lll}{{\boldsymbol{\left[\right.}\!\!\boldsymbol{\left[\right.}}}
\newcommand{\rrr}{{\boldsymbol{\left.\right]}\!\!\boldsymbol{\left.\right]}}}

\renewcommand{\o}{\otimes }

\renewcommand{\k}{{\C}}
\newcommand{\Z}{{\mathbb{Z}}}

\newcommand{\LL}{{\mathbb{L}}}
\newcommand{\LS}{{\mathbb{L}^\sharp_\dd}}
\newcommand{\MM}{{\mathbb{M}}}
\newcommand{\Int}{{\mathbf{I}}}

\newcommand{\ze}{\zeta}
\newcommand{\UU}{{U}}
\newcommand{\VV}{{\mathbb{V}}}
\newcommand{\EE}{{\mathbf{E}}}

\newcommand{\K}{{K}}
\newcommand{\KK}{K^{\operatorname{cyc}}}
\newcommand{\BK}{{\mathbf K}}
\newcommand{\even}{{\operatorname{even}}}
\newcommand{\odd}{{\operatorname{odd}}}

\newcommand{\ben}{\begin{enumerate}}
\newcommand{\een}{\end{enumerate}}

\newcommand{\super}{_{\text{super}}}
\newcommand{\anti}{_{\text{anticyc}}}

\newcommand{\ev}{\operatorname{ev}}

\newcommand{\CC}{{\mathbb{C}}}
\newcommand{\C}{{\mathbb{C}}}
\newcommand{\N}{{\mathbb{N}}}
\newcommand{\ZZ}{{\mathbb{Z}}}

\newcommand{\pa}{\partial }
\newcommand{\lo}{^\circ }
\renewcommand{\SS}{{\mathscr S}}
\theoremstyle{plain}
\newcommand{\TR}{{\mathbf{Tr}}_{\dd}}
\newcommand{\TS}{{\mathbf {Tr}}\lo }
\newcommand{\dd}{{\mathbf{d}}}

\newcommand{\rda}{{\Rep_\dd{A}}}
\newcommand{\QQ}{{\overline{Q}}}
\def\hp{\hphantom{x}}
\newcommand{\ww}{{\mathbf{w}}}

\newcommand{\bp}{{\mathbf{p}}}
\newcommand{\bq}{{\mathbf{q}}}
\newcommand{\bg}{{\mathbf{g}}}
\newcommand{\WW}{{\mathbf{W}}}
\newcommand{\PP}{{\mathbf{P}}}
\newcommand{\BQ}{{\mathbf{Q}}}
\newcommand{\FM}{{\mathfrak M}}

\newcommand{\bone}{{\boldsymbol{1}}}
\begin{document}

\title{\qquad Noncommutative complete intersections \newline
and matrix integrals}

\author{Pavel Etingof}
\address{Department of Mathematics, Massachusetts Institute of Technology,
Cambridge, MA 02139, USA}
\email{etingof@math.mit.edu}

\author{Victor Ginzburg}
\address{Department of Mathematics, University of Chicago,
Chicago, IL 60637, USA}
\email{ginzburg@math.uchicago.edu}

\maketitle
\centerline{\textbf{\em To Bob MacPherson on the occasion of his 60th
birthday}}

\begin{abstract} We introduce a class 
of noncommutatative algebras called {\em  representation
complete intersections} (RCI). A graded
 associative algebra $A$ is said to be RCI  provided
there exist arbitrarily large positive integers $n$ such that
 the scheme $\Rep_nA$,  of $n$-dimensional representations
of $A$,  is a complete intersection.
 We discuss  examples of RCI algebras, including
those arising from quivers.

There is  another interesting class of
 associative algebras called {\em   noncommutative  complete
intersections} (NCCI). We
prove  that any graded RCI algebra  is NCCI. We also
 obtain explicit formulas for
the Hilbert series of each nonvanishing cyclic and Hochschild homology group
of an RCI algebra.
The proof involves a  noncommutative
 {\em   cyclic Koszul complex}, $\KK_\idot{A},$
and   a matrix integral similar to the one arising
in  quiver gauge theory.
 \end{abstract}
\vskip 3pt

\centerline{\sf Table of Contents}
\vskip -1mm

$\hspace{30mm}$ {\footnotesize \parbox[t]{115mm}{\hp${}_{}$\!
\hp\!1.{ $\;\,$} {\tt Introduction} \newline
\hp2.{ $\;\,$} {\tt Representation functor and matrix integrals}\newline
\hp3.{ $\;\,$} {\tt Noncommutative complete intersections}\newline
\hp4.{ $\;\,$} {\tt The trace map}\newline
\hp5.{ $\;\,$} {\tt Additional results and examples}\newline
\hp6.{ $\;\,$} {\tt Preprojective algebras and quiver varieties}
}
}
\section{{Introduction}}
\subsec{}
In this paper we work with finitely generated unital associative 
algebras over $\mathbb C$. Such an algebra $F$ is said to be
{\em smooth} if $\Ker[F\o F\to F]$, the kernel of
the multiplication map, is a projective $F$-bimodule.

Fix  a smooth algebra $F$, for instance a free algebra
$F={\mathbb C}\langle x_1,\ldots, x_d\rangle.$ Let $J$ be a two-sided finitely generated
ideal in $F.$

\begin{definition} The algebra $A=F/J$ (or the pair $J\sset F$)
is called a {\em noncommutative
complete intersection (NCCI)}  if $J/J^2$ is a projective
$A$-bimodule. This is equivalent to the condition that $A$ has
a projective $A$-bimodule resolution of length $\leq 2,$ cf. Theorem~\ref{ncci}.
\end{definition}

NCCI algebras have been considered, under various different names,
by a number of authors, see e.g. \cite{AH}, \cite{An}, \cite{GS}, \cite{Go}, \cite{Pi}.
In the present paper, we are interested in the case of {\em graded}
NCCI algebras. Specifically, fix a finite set  $I$
and write $R$ for the algebra of functions $I\to\CC,$
with pointwise multiplication. Let
$Q$ be a finite  quiver with vertex set $I$, and
assign  to each  edge $a\in Q$ an arbitrary positive
grade degree $d(a)>0$.  This
makes  the path algebra of  $Q$ a  free graded  $R$-algebra;
moreover, any  free graded  $R$-algebra $F$ may be obtained
in this way.

We study
algebras of the form
$A=F/J,$ where $J$ is a graded
two-sided ideal in the path algebra
$F$ such that $J/J^2$ is a projective finitely generated
$A$-bimodule.
Thus, $A=\bigoplus_{r\geq 0}A[r]$ is
 a graded   algebra such that  $A[0]=R$.

The grading on  $A$ induces
a natural grading on each
Hochschild  homology group ${HH}_k(A)$.
Our main results  (Theorem \ref{maint} and Theorem \ref{converse})
give
 explicit   formulas
for the Hilbert series of each Hochschild homology group of $A$
 for a class of algebras $A$ that we call
{\em asymptotic
representation complete intersections} (see below). 

For $k=0$, for instance, Hochschild homology reduces to
the commutator  quotient space,  ${HH}_0(A)=A/[A,A],$
where $[A,A]$ denotes the $\k$-linear subspace
of $A$ spanned by the
commutators $ab-ba,\,a,b\in A$. 
The grading on $A$  clearly descends to the commutator  quotient.
Our result yields a formula,
in terms of an infinite product of determinants,
for the  Hilbert series of the graded algebra
\beq{OA}
\O(A)=\sym_\C\left({A_+}\big/{[A,A]}\right).
\eeq
Here,  $A_+=\bigoplus_{r> 0}A[r]$,
so  $A_+/[A,A]$ is  the positive degree
part of $A/[A,A]$.
We remark that the symmetric algebra  in \eqref{OA}
is equipped with 
the {\em total} grading 
that comes from the  grading on  $A_+/[A,A]$
(this grading should not be confused with the standard grading 
$\sym =\bigoplus_{k\geq 0} \sym^k,$ on the Symmetric algebra of a vector space).

\subsec{} The main idea of our approach is based on the notion of
{\em representation functor} and
 may be outlined as follows. 

Let $A$ be any graded finitely presented 
algebra. 
For each positive  dimension  $\dd$,
one has an affine scheme $\Rep_\dd{A}$
of $\mathbf d$-dimensional representations of $A$, cf. \S\ref {basic}.
Heuristically, one expects that the
collection  of representation schemes  $\Rep_\dd{A}$
provides, as $\dd\to\infty,$  a `good approximation' of
the algebra $A$ (this is not true for nongraded algebras which
may have no finite dimensional representations at all).
 In particular, one should be able to read off
a lot of information about $A$ from some geometric
information about the corresponding representation schemes.

In more detail, for each $\dd$,
the general linear group  $GL_\dd$ acts naturally
on the scheme $\Rep_\dd{A}$ by  base change transformations 
of the underlying
vector space of the
$\mathbf d$-dimensional representation.
We consider the algebra $\C[\Rep_\dd{A}]^{GL_\dd},$
of $GL_\dd$-invariant 
polynomial functions on  $\Rep_\dd{A}$.
We show that, in some  sense,
the commutative algebra $\O(A)$, in \eqref{OA},
is a limit of the algebras $\C[\Rep_\dd{A}]^{GL_\dd}$ as $\mathbf d$ goes to
infinity. 
A precise version of this statement will be proved in \S\ref{sec3}.
The statement implies, in particular, that the Hilbert series of
 $\O(A)$ may be obtained as a limit of 
  Hilbert series of the algebras  $\C[\Rep_\dd{A}]^{GL_\dd}$ as $\mathbf d$ goes to
infinity. 

Observe further that the presentation $A=F/J$
makes   $\Rep_\dd A$ a closed subscheme
in $\Rep_\dd F$, the latter being a vector space. We 
 introduce a noncommutative  {\em cyclic Koszul complex},
$\KK_\idot{A},$ such that each of the
 functors $\Rep_\dd,\,\dd\in\Z^I,$
sends $\KK_\idot{A}$ to the $GL_\dd$-fixed part of the standard
Koszul complex of the subscheme
$\Rep_\dd A\sset \Rep_\dd F$. Moreover, we prove that
the two complexes become `asymptotically isomorphic'
as $\dd\to\infty.$

To be able to exploit usefully the `asymptotic isomorphisms'
above we introduce the notion of
 {\em  representation   complete intersection} (RCI) algebra,
and a weaker (but more useful)  notion of
{\em asymptotic} RCI algebra.
The algebra $A=F/J$ is called RCI provided
there exist arbitrarily large positive dimensions $\dd
$ such that
the scheme $\Rep_\dd A$
 is a complete intersection in  $\Rep_\dd F$.
We will prove:
\beq{equiv}
\text{RCI}\en\;\Rightarrow\en\;
\text{asymptotic RCI}\en\;\Leftrightarrow
\en\;
\text{NCCI with `expected' }\; {HH}_2(A)\en\;\Rightarrow
\en\;\text{NCCI}.
\eeq

It turns out that, for an asymptotic  RCI algebra $A$,
there is enough  asymptotic information
about Koszul complexes of the schemes
 $\Rep_\dd A$ to deduce formulas for the Hilbert series
of each Hochschild homology group of $A$.

\subsec{Quiver algebras.}\label{uninak} 
Our results apply, in particular,
in an interesting special case arising
from quivers. 
Specifically, let $Q$
be  a (nonoriented)  graph
with vertex set $I$. Let $\QQ$ be the quiver
obtained by doubling  edges of $Q$ so that 
each nonoriented edge of $Q$ gives rise to a pair, $a,a^*,$ of
oriented edges 
of $\QQ$ equipped with opposite orientations.
Write $\bc=\|\bc_{ij}\|$ for the adjacency  matrix
of  $\QQ$, thus $\bc_{ij}$ denotes the number of
edges $i\to j$ in $\QQ$.

We take $F=F(Q)$ to be the path algebra of $\QQ$.
Let
$$\Pi= F(Q)\Big/\big(
\sum\nolimits_{a\in Q}
 [a,a^*]\big)
$$
be the {\em preprojective algebra} associated to $Q$,
a quotient of $F(Q)$ by the two-sided ideal  generated by the
element
$\sum_{a\in Q}
 (a\cdot a^*-a^*\cdot a)\in F(Q).$ 
 The grading on the path algebra induces
natural gradings on $\Pi$ and on $\Pi/[\Pi,\Pi].$

Write $\bone$ for the identity $I\times I$-matrix.
Thus, $\bone-t\cd\bc+t^2\cd\bone$ is an $I\times I$-matrix
with entries in
$\ZZ[t].$ This matrix specializes at $t=1$ to the 
 {\em Cartan matrix} of  $Q$. 
The closely related matrix $\frac{1}{t}(\bone-t\cd\bc+t^2\cd\bone)=
(t+t\inv)\cd\bone-\bc$, called $t$-analogue
of the  Cartan matrix, has
been considered by Kostant and Lusztig, see \cite[\S3.12-3.13]{Lu}.

Applying our general theorem to the algebra $\Pi$ yields the following result.

\begin{theorem}\label{nak} Let $Q$ be neither Dynkin nor extended 
Dynkin quiver.
Then, $\Pi=\Pi(Q)$ is an RCI algebra, and the Hilbert series of $\Pi$,
resp. of
$\Pi/[\Pi, \Pi]$, is determined from the formula
$$
h(\Pi;t)=(\bone-t\cd\bc+t^2\cd\bone)\inv,\quad\text{resp.}\quad
h\big(\O(\Pi);t\big)=\frac{1}{1-t^2}
\prod_{r\ge 1}\frac{1}{\det (\bone-t^r\cd\bc+t^{2r}\cd\bone)}.
$$
\end{theorem}
Here, the formula for $h(\Pi,t)$ has been known for some time, see
e.g. \cite{MOV} and references therein,
while the formula for $h\big(\O(\Pi); t\big)$ is new.

Recall that, according  to Van den Bergh
\cite{VdB}, there are natural isomorphisms
${HH}_\idot(\Pi)$ $\cong {HH}^{2-\hdot}(\Pi)$ relating Hochschild
homology and {\em co}homology.
These isomorphisms and the
 exact sequence \eqref{exact},
see Theorem \ref{maint} below,
 yield 
formulas for Hilbert series of all Hochschild {\em co}homology
groups of the algebra $\Pi$.
In particular, this gives an alternate proof of
the result that, for any quiver of  neither Dynkin nor extended 
Dynkin type, the center of the algebra $\Pi$ reduces to scalars,
cf. \cite[Proposition 8.2.2]{CBEG}.

Both  formulas of  Theorem \ref{nak} fail for  $A,D,E$ Dynkin  quivers.

Now, let $Q$ be an  {\em extended} 
Dynkin quiver with extending vertex $o\in I$.
Then, the  formula for $h(\Pi;t)$ given in  Theorem \ref{nak} still holds,
but the formula for $h\big(\O(\Pi);t\big)$ has to be modified,
see Corollary \ref{b} in ~\S\ref{extended},
and  \cite{Su}.
In particular, that modified formula implies
a curious identity for the
adjacency matrix of $\QQ$
which we were unable to find in the literature.

The identity involves  {\em Chebyshev polynomials}
 of the first kind. Specifically, put $\varphi_0:=1$ and let
$\varphi_k\in\C[x],\,k=1,2,\ldots,$ be
 a  sequence 
of polynomials in an indeterminate $x$
 defined as the coefficients in
 the  expansion of the generating function
\beq{varphi} 
\frac{2-tx}{1-tx+t^2}=1+\sum_{k=1}^\infty
\varphi_k(x)t^k,
\quad\text{thus}\quad
\varphi_k(2\cos{z}):=2\cos(k z),\en\forall k\geq 1.
\eeq

For each $k\geq0,$ we may plug the  adjacency matrix $\bc$,  of $\QQ$,
into the polynomial $\varphi_k$.
 This way, we obtain
 an integer valued $I\times I$-matrix
$\varphi_k(\bc)=\|\varphi_k(\bc)_{ij}\|$. 
The matrix entries corresponding to
the extended vertex $i=j=o$ form 
 a sequence of integers $\varphi_k(\bc)_{oo},\,k=0,1,\ldots.$
Our identity reads
\beq{curious}
\prod_{r\geq 1}\det (\bone-t^r\cd \bc+t^{2r}\cd \bone)=
\prod_{k\geq 1}(1-t^k)^{\varphi_k(\bc)_{oo}}.
\eeq

This identity will be proved in \S\ref{extended} by applying the Euler-Poincar\'e
principle to $\KK_\idot\Pi$, the noncommutative  cyclic Koszul complex of the
preprojective algebra $\Pi$.
One can also verify the identity
 case by case,
 for each extended Dynkin quiver of $A,D,E$ types separately,
 see \S\ref{extended} for details.

\subsec{Layout of the paper.} In Section 2, we  formulate our 
first important result (Theorem \ref{maint})
 and carry out crucial
 matrix integral calculations.
The notion of RCI algebra is introduced in
\S\ref{setting} and the notion of asymptotic
RCI algebra  is introduced in
\S\ref{asymptotic}.
We use the standard  Koszul
complex associated with
 a complete intersection to  express the Hilbert series  of the
 coordinate ring   $\C[\Rep_\dd{A}]$ in terms
of a matrix integral,  an integral over the unitary group $\UU_\dd$.
Then, in \S\ref{asymptotic}, we apply some results
from the theory of random matrices to 
find the asymptotics of our  matrix integral
as $\dd\to\infty$.  For an  asymptotic
RCI algebra, that enables us to obtain an asymptotic  formula for 
 the   Hilbert series  of the algebras
$\C[\Rep_\dd{A}]$ as $\dd$ goes to infinity.

In \S\ref{sec2} we discuss NCCI algebras and related
 homological algebra. Similarly to the 
familiar case of commutative algebras,  to any noncommutative algebra
 $A$
given by generators and relations one can associate canonically
a free  DG algebra $K_\idot{A},$ called the
(noncommutative)
Koszul complex of $A$. The noncommutative
Koszul complex is defined in \S\ref{kos}.
 It turns out that $A$ is an  NCCI algebra
if and only if  the DG algebra
 $K_\idot{A}$  is quasi-isomorphic to $A$.
Main results about NCCI algebras are gathered
in Theorems \ref{ncci} and  \ref{Anick}.

In \S\ref{CCI}, we 
  introduce the cyclic 
Koszul complex $\KK_\idot{A}$ that computes cyclic homology
of the NCCI  algebra $A$. We prove that these vanish
in all degrees $>1$. We formulate
the main result of the paper, Theorem \ref{converse},
which claims, among other things, the equivalence in the
middle of \eqref{equiv}.

The proofs of the two main theorems are completed in
\S\ref{sec3}. In that section, we compare the
algebras $A$ and $\C[\Rep_\dd{A}]$ via  natural evaluation
maps. In particular, we construct  an evaluation morphism
from
the  cyclic 
Koszul complex  $\KK_\idot{A}$ to the ordinary Koszul complex  
of the representation scheme $\Rep_\dd{A},$
that commutes with the respective Koszul differentials.
We prove that, in the $\dd\to\infty$
limit,   evaluation morphisms become
isomorphisms. Thus, homology vanishing 
for the noncommutative  Koszul complex   $\KK_\idot{A}$ is equivalent
to an asymptotic  homology vanishing for the ordinary
Koszul complexes of  representation schemes $\Rep_\dd{A}.$
Furthermore the   Hilbert series  of $\KK_\idot{A}$
may be expressed in terms of  Hilbert series  
of the corresponding Koszul complexes of  representation schemes.

In \S5, we consider various examples of RCI and asymptotic RCI algebras.
We also provide examples of non RCI algebras which, nevertheless,
are asymptotic RCI, see  \S\ref{asex}. Specifically,
 we show that  if 
the relations of $A$ are pairwise non-overlapping monomials,
then $A$ is an asymptotic RCI but not necessarily RCI.
We  use a deformation argument to deduce a similar result
for algebras with
 generic  homogeneous relations of degrees equal to the
lengths of these non-overlapping monomials. 
This way, one obtains many nontrivial examples of  asymptotic RCI
algebras with explicit
infinite product formulas for the corresponding Hilbert
series. 

In \S6, we show
 that, for   connected quivers of neither Dynkin nor extended
Dynkin type, the corresponding preprojective algebra is RCI.
Also,  partial preprojective algebras of any connected
quiver (in the sense of \cite{EE}),  are RCI algebras.
In \S\ref{mdw}, we give an interesting asymptotic formula for Hilbert series of
Nakajima quiver varieties. Extended Dynkin quivers
are discussed in \S\ref{extended}, where we prove the identity in
\eqref{curious}.

\begin{remark}\label{not} In the present paper, the ground field is
 the field $\C$ of complex numbers.
All the results of the paper
can be routinely
generalized to the setting of an arbitrary
algebraically closed ground field of characteristic zero.
Characteristic zero assumption is essential
for the geometric results involving representation schemes.
It is also essential for the cyclic/Hochschild
 homology computation in Lemma \ref{hclem}
and Proposition \ref{HC1},
which involves Karoubi-de Rham complex.
However, the rest of \S3 applies verbatim in the
case of ground fields of arbitrary characteristic.
\end{remark}
\vskip -6pt

\subsec{Ackmowledgments.}{{\footnotesize{The authors are grateful 
to M. Artin, L. Hesselholt, G. Lusztig, H. Nakajima,  D. Piontkovsky, E. Rains,
and R. Stanley 
for useful comments and discussions. We also thank E. Rains for explaining to us the results 
of Diaconis and Shahshahani on matrix integrals, 
and for providing a MAGMA code that we used at all stages of this
work. The work of P.E. and V.G.
 was  partially supported by the NSF grants
 DMS-0504847 and DMS-0303465,
 respectively, and  by the CRDF grant RM1-2545-MO-03.}}}

\section{{ {Representation functor and matrix integrals}}}\label{secbas}
\subsec{Basic definitions.}\label{basic}
Throughout, we fix  a finite set $I$ and let
$R$ be the algebra of functions $I\to\k,$
with pointwise multiplication. Thus, $R$ is a commutative
 semisimple $\k$-algebra.
Simple 
$R$-modules are 1-dimensional and are parametrized by 
elements of the set $I$; an $R$-module is the same thing as
an $I$-graded vector space, $M=\oplus_{i\in I}M_i$.
If $M$ is finite-dimensional and $d_i=\dim M_i$, we
say that $\dd=\{d_i\}_{i\in I}$, an  $I$-tuple of nonnegative integers,
is the
{\em dimension vector} for $M$.

Similarly, $R$-bimodules are $I\times I$-graded vector spaces, 
$M=\oplus_{i,j\in I}M_{ij}$. 

Let  $M=\bigoplus_{r\geq 0} M[r]$ be
 a $\mathbb Z_+$-graded $R$-bimodule
such that  each graded component $M[r]$ is a finite dimensional
 $R$-bimodule, that is we have
 $M[r]=\bigoplus_{i,j\in I} M_{ij}[r],$
where $\dim  M_{ij}[r]<\infty.$
The Hilbert  series of $M$ is a $\ZZ[[t]]$-valued
$ I\times I$-matrix
defined as
$$h(M)=\|h_{ij}(M; t)\|,\quad
h_{ij}(M; t)
:=\sum\nolimits_{r\ge 0}\,\dim M_{ij}[r]\cd t^r.
$$
It is clear that  $h({M\otimes_R M'})=h(M)\cdot h(M')$
for any  $\mathbb Z_+$-graded 
$R$-bimodules  $M$ and~$M'$.

Given a dimension vector $\mathbf d$, let
$\C^{\dd}:=\bigoplus_{i\in I}\C^{d_i}$.
This is an $I$-graded 
vector space, and we consider an associative algebra
$$\EE_\dd :=\Hom_\k(\C^{\dd}, \C^{\dd})=\bigoplus\nolimits_{i,j\in I}\EE_{ij},
\quad\text{where}\quad
\EE_{ij}=\Hom_\k(\C^{d_i},\C^{d_j}).$$
The $I$-grading makes the vector space $\C^{\dd}$  a left $R$-module,
and the action map for this module gives a natural algebra homomorphism
$R\to \bigoplus_{i\in I}\EE_{ii}\sset \EE_\dd .$

\begin{notation}
We use unadorned symbols like
$\sym, \La, \Hom, \o, $ etc.,  for  symmetric algebra,
exterior algebra, $\Hom$-space,
tensor product, etc., all
taken over  $\C$ (not over $R$).

Given 
a finite dimensional  $R$-bimodule
$M$, we reserve the corresponding blackboard font notation
for the vector space 
\beq{bold}
\MM_\dd:=\Hom_{\bimod{R}}(M,\EE_\dd ).
\eeq 
The exceptions to this convention are:
$\N,\Z_+,\Z,$ and $\C$, the sets  of natural numbers, nonnegative integers,
integers, and
complex numbers,
respectively.\hfill$\lozenge$
\end{notation}

An $\N$-grading on  an $R$-bimodule
$M$ gives rise to a linear $\C^\times$-action
on the  vector space $\MM_\dd$. This action
contracts the vector space to its origin,
 the only
$\C^\times$-fixed point.

Let $F$ be an $R$-algebra, that is, a
{\em finitely presented} associative unital $\k$-algebra
equipped with an algebra imbedding $R\into F$.
Given 
a dimension vector $\mathbf d$, let
$\Rep_{\mathbf d}F$
be the set of all algebra maps
$F\to \EE_\dd :=\Hom_\k(\C^{\dd}, \C^{\dd})$ which restrict to the
natural map $R\to \EE_\dd .$ The
set $\Rep_{\mathbf d}F$ has the natural structure
of a (not necessarily reduced) affine scheme of finite type,
 called the  {\em representation scheme} of 
$\dd$-dimensional representations of the algebra $F$.
Let $\k[\Rep_\dd{F}]$ denote the coordinate ring of
this scheme.

\subsec{Main result.}\label{setting}
Let  $T_RV=R \bigoplus V \bigoplus V\otimes_R V\bigoplus\ldots,$ denote the tensor algebra of an $R$-bimodule
$V,$
and write $T_R^+V=V \bigoplus V\otimes_R V\bigoplus\ldots,$ for  the augmentation ideal of $T_RV$.

\begin{definition} A pair of  finite dimensional $\mathbb N$-graded
$R$-bimodules
$V$ and $L$, together with a graded $R$-bimodule  imbedding
 $j: L\into T_R^+V,$ will be referred to as  $(V,L)$-{\em datum}.
\end{definition}

Fix a  $(V,L)$-datum, where
 $V=\bigoplus_{r\geq 1} V[r],$ and
$L=\bigoplus_{r\geq 1} L[r].$
We define 
\beq{lo}L\lo:=j\inv\big([T_RV,T_RV]\big)\sset L.
\eeq
Thus, $L\lo=\bigoplus_{r\geq 1} L\lo[r]$ is a graded
subspace in $L$. We put $m_r=\dim L\lo[r],$
and let 
\beq{lambda}
\lambda(L\lo) :=\prod\nolimits_{r>0}\,
(1-t^r)^{m_r}\in \ZZ[t].
\eeq

We also introduce the following generating $\zeta$-function
\beq{zeta}
\ze(V,L):=
\prod_{s=1}^\infty \frac{1}{\det\big(\bone- h(V; t^s)+h(L; t^s)\big)} \in
\ZZ[[t]].
\eeq

The  gradings on $V,L,$ and $L\lo,$ give rise to  {\em total} gradings on 
various  objects like exterior algebra $\La{L\lo}$,
tensor algebra $T_RV$, etc.
These  total gradings  will be called
{\em weight} gradings. For the corresponding bigraded algebras, we write 
e.g., $\La L\lo=\bigoplus_{p,r_{_{}}\geq 0}$
$\dis(\La^p_{_{}} L\lo)[r].$
With this notation, for the polynomial in \eqref{lambda}, one has
$$\lambda(L\lo)=\sum\nolimits_{p,r}(-1)^p\cdot t^{r^{}}\cdot \dim ((\La^p L\lo)[r]).$$

Associated with a $(V,L)$-datum, is the
two-sided ideal 
$J:=(L)\sset T_RV$ generated by $j(L)$,
and the corresponding quotient
algebra  $A:=T_RV/J=T_RV/(L).$ 
The weight grading on $T_RV$ makes $J$ 
 a graded ideal and $A$ a graded algebra,
$A=\bigoplus_{r\geq 0}A[r]$. For each $r\geq 0,$
the homogeneous component
$A[r]$ is a finite dimensional $R$-bimodule;
we have $A[0]=R$.

Given an integer $N$ and a dimension vector $\dd$, we write $\dd \succ N$ if
$d_i > N$ for all $i\in I.$

\begin{definition}\label{rci_def} A $(V,L)$-datum is called
{\em representation complete intersection (RCI)}
if for any integer $N$ there exists a
dimension vector $\mathbf d\succ N$ such
that, cf.  \eqref{bold}:
\beq{dim}\dim\Rep_\dd A=\dim\VV_\dd -\dim\LL_\dd +\dim L\lo.
\eeq
\end{definition}

The geometric meaning of equation \eqref{dim} will be explained in
\S\ref{hilb_rci} below.

Abusing terminology, we often say 
that the algebra $A=T_RV/(L)$  arising from
the  $(V,L)$-datum is an RCI algebra.

One of our main results reads

\begin{theorem}\label{maint} Let $(V,L)$ be an RCI datum and
$A=T_RV/(L)$. Then, we have ${HH}_k(A)=0$ for all $k>2$ (i.e., $A$ is
a NCCI algebra, cf. \S\ref{NCCI}), and ${HH}_2(A)\cong L\lo$. 

Furthermore,  there is a natural
exact sequence
\beq{exact}
0\too R\too {HH}_0(A)\too {HH}_1(A)\too {HH}_2(A)\too 0.
\eeq

The   Hilbert
series of 
$A$, resp. of ${HH}_0(A)=A/[A,A],$ is determined by the formula
\beq{AA}
{\mathsf{(i)}}\quad h(A)=\big(\bone-h(V)+h(L)\big)\inv,\quad\text{resp.},\quad
{\mathsf{(ii)}} \quad h(\O(A))=\ze(V,L)/\lambda(L\lo) .
\eeq
%
\end{theorem}

\begin{remark} One may call
${\mathbf{p}}(V,L;t):=\bone-
 h(V; t)+h(L; t)$, a $\ZZ[t]$-valued $I\times I$-matrix,
the
{\em Cartan polynomial} of the $(V,L)$-datum. Thus,  formulas
\eqref{AA} read
$$h(A)={\mathbf{p}}(V,L;t)\inv,\quad\text{resp.},\quad
h(\O(A))=\frac{1}{\lambda(L\lo)}\cd\prod_{s=1}^\infty 
\frac{1}{\det{\mathbf{p}}(V,L;t^s)}.
$$
\end{remark}

\begin{remark} It follows from 
the isomorphism $L\lo\cong {HH}_2(A)$ and from \eqref{exact}
that one has $h({HH}_1(A))$
$=h\big(A_+\big/[A,A]\big)+ h( L\lo).$
Therefore, the Theorem implies also  that the Hilbert
series of ${HH}_1(A)$ may be found from the formula
$$h(\sym {HH}_1(A))=\ze(V,L)/\lambda(L\lo)^2.$$
\end{remark}

The proof of Theorem \ref{maint} will be completed in \S\ref{sec3}.
We provide, in effect, two different proofs of the Theorem.
The first proof is based on a matrix integral calculation,
to be carried out in \S\S\ref{hilb_rci}-\ref{asymptotic} below.
The second, purely algebraic proof, is based on the
formula for $h(\O(A))$ in the special case of
a {\em free} $R$-algebra $A$. Such a formula is known,
it can be obtained
by a combinatorial argument that involves counting 
cyclic paths on a graph, cf. \cite{St}. The general
case of  Theorem \ref{maint} for an arbitrary RCI algebra
$A$ can then be deduced from the special case
of free $R$-algebras by applying  the 
formula proved in the special case to the noncommutative
Koszul complex $K_\idot{A}$, which is free as
an $R$-algebra, cf. ~\S\ref {kos}.

\subsec{Hilbert series of representation schemes.}\label{hilb_rci}
Fix 
 a dimension vector $\mathbf d$ and a $(V,L)$-datum.
Clearly,  one has a natural isomorphism
$\Rep_\dd(T_RV)\cong \Hom_{\bimod{R}}(V,\EE_\dd)=\VV_\dd,$
cf. \eqref{bold}.
 
Now,
let $A=T_RV/(L)$.
The algebra projection $T_RV\onto A$
 induces a natural closed
imbedding of representation schemes $\Rep_\dd A\into\Rep_\dd(T_RV).$
Dually,
the  $R$-bimodule imbedding $j: L\into  T_RV$ induces a natural
restriction morphism 
\beq{j^*}
j^*_\dd:\ \VV_\dd 
=\Hom_{\bimod{R}}(V,\EE_\dd )\cong\Rep_\dd(T_RV) 
\too\LL_\dd.
\eeq
It is clear that the scheme $\Rep_\dd A$ is 
 the scheme-theoretic zero fiber of the morphism
$j^*_\dd$.

In general, the map
$j^*_\dd$ is {\em not} surjective.
Indeed,  for any algebra homomorphism
$\rho: T_RV\to \EE_\dd$ and any $x\in [T_RV,T_RV],$
we have $\Tr\rho(x)=0.$ In particular, 
the composite map 
$$L\lo\stackrel{j}\into [T_RV,T_RV]\into T_RV\stackrel{\rho}\too
\EE_\dd \stackrel{\Tr}\too\C
$$
must vanish. Therefore, the image of the map
$j^*_\dd$ in \eqref{j^*} is contained in a proper
subspace
$\LS:=\Ker\psi\sset\LL_\dd$, the
kernel of the following composite map
\beq{LS}
\LS:=\Ker\big[\psi:
\xymatrix{
\LL_\dd\;
\ar@{->>}[r]^<>(0.5){\text{res}^L_{L\lo}}&\;\Hom(L\lo,\EE_\dd )\;
\ar@{->>}[r]^<>(0.5){\Tr|_{\EE_\dd}}&\;\Hom(L\lo,\C)=(L\lo)^*
}\big].
\eeq

Write $p_\dd:  \Rep_{\mathbf d}(T_RV) \to \LS$
for the resulting map
and let $o$ denote the origin
of the vector space $\LS$. This way, we get the
following commutative diagram

\beq{flat}
\xymatrix{
\Rep_{\mathbf d} A\ar@{^{(}->}[rr]\ar[d]&&\VV_\dd =\Rep_{\mathbf
d}(T_RV)
 \ar[d]^<>(0.5){p_\dd}\ar[drr]^<>(0.5){j^*_\dd}&\\
\{o\}\ar@{^{(}->}[rr]&&\LS\ar@{^{(}->}[rr]&&
\LL_\dd
}
\eeq

The  left square in the diagram is cartesian,
so  $\Rep_\dd A$ is the scheme-theoretic zero fiber
of the morphism $p_\dd$. The following result
provides a geometric meaning for the notion of
representation complete intersection datum.

\begin{lemma}\label{meaning} The dimension equality \eqref{dim} holds
 for
a dimension
vector $\dd$  if and only if the
map $p_\dd$, in
diagram \eqref{diag}, is {\em flat}.
 \end{lemma}

To prove the Lemma, we observe  that the
gradings on $V$ and $L$ make each of the 
spaces in diagram \eqref{flat} a   $\C^\times$-variety, such that
all maps in the diagram become $\C^\times$-equivariant
morphisms. Moreover,  the $\C^\times$-action on $\LS $ being a contraction,
we deduce that the dimension of any fiber of the map $p_\dd$
is less than or equal to the dimension of the zero fiber,
that is to $\dim\Rep_\dd A.$
Therefore,
 the morphism $p_\dd$  is flat if and only if one has
$\dim\Rep_\dd A=\dim \VV_\dd-\dim \LS.$

Observe next that the map $\psi$ in  \eqref{LS}
is surjective. Hence,
one has an exact sequence
\beq{LS_exact}
0\too \LS\too \LL_\dd\stackrel{\psi}\too (L\lo)^*\too 0.
\eeq
We deduce:
$\dim \LS=\dim \LL_\dd-\dim  L\lo$, hence,
$\dim \VV_\dd-\dim \LS=\dim \VV_\dd-\dim \LL_\dd+\dim  L\lo$, 
and the Lemma
follows.
\qed
\vskip 3pt

Given $i\in I$ and an element $g_i\in GL(\C^{d_i})$,
we write $g_i^\vee=(g_i^\top)\inv: (\C^{d_i})^*\to (\C^{d_i})^*,$
where $g_i^\top$ denotes the dual endomorphism of the dual vector space.
For each pair of vertices
$i,j\in I,$ 
we consider  endomorphisms of  the vector space
$(\C^{d_i})^*\o\C^{d_j}=\Hom(\C^{d_i},\C^{d_j})=\EE_{ij}$
of the form
$g_i^\vee\o g_j$.
Thus, given $g_i\in GL(\C^{d_i})$ and $g_j\in GL(\C^{d_j})$, for any integer $r\geq 0$,
there is a well defined  polynomial
$$
\det(\bone_{\EE_{ij}}-t^r\cdot g_i^\vee\o g_j)\in\C[t].$$

We put
$G_\dd=\prod_{i\in I} GL(\C^{d_i}).$ 
This is a reductive group that acts naturally  on
the vector space $\C^{\dd}=\bigoplus_{i\in I}\C^{d_i},$ hence, 
also on the algebra $\EE_\dd,$
by conjugation.

Next, for each $d,$ we fix a maximal compact subgroup $U(d)\sset
GL(\C^d).$ 
The group $\UU_\dd:=
\prod_{i\in I}U(d_i)$
is a maximal compact subgroup of $G_\dd$.
Let $dg$ be the Haar measure on $\UU_\dd $ normalized so
that the total volume of  $\UU_\dd $ equals 1.

We introduce
 the following matrix integral
\beq{Int}
\Int(V,L,\dd):=\int_{\UU_\dd}\;\prod_{r>0}\prod_{i,j\in I}\;
\frac{dg}{\det(\bone_{\EE_{ij}}-t^r\cdot g_i^\vee\o g_j)^{\bc_{ij}^r(V,L)}}.
\eeq
Here,  $\bc_{ij}^r(V,L):=\dim V_{ij}[r]-\dim L_{ij}[r]$,
and the above expression is viewed as a formal power series 
$\Int(V,L,\dd)=\sum_{m\geq 0}\Int_m(V,L,\dd)\cd t^m\in \C[[t]]$ given by the 
 power series expansion of the RHS of \eqref{Int}.

Conjugation-action of the group $G_\dd$ on $\EE_\dd$ induces  linear $G_\dd$-actions
on the vector spaces $\VV_\dd,\LL_\dd,$ etc. 
Furthermore, all the maps in
diagram \eqref{diag} are $G_\dd$-equivariant morphisms.
In particular, the scheme
$\Rep_\dd A$ is a $G_\dd$-stable subscheme
of $\VV_\dd$,
so the group $G_\dd$ acts on the coordinate ring $\k[\Rep_\dd{A}]$
by algebra automorphisms. 

We also have the
 $\C^\times$-action on the scheme $\Rep_\dd{A}\sset\VV_\dd$
that commutes with the $G_\dd$-action.
The induced $\C^\times$-action on
the subalgebra of $G_\dd$-invariants gives rise to a grading
\beq{weight}
\k[\Rep_\dd{A}]^{G_\dd}=\bigoplus\nolimits_{r\geq 0}\,
\k[\Rep_\dd{A}]^{G_\dd}[r],
\eeq
 to be referred to as  {\em weight grading}.

\begin{proposition}\label{intform} Let $A=T_RV$ be the graded
algebra associated to a $(V,L)$-datum. Assume that the map $p_\dd$ in
\eqref{flat}
is flat, so the scheme $\Rep_\dd{A}$
is a complete intersection in $\VV_\dd $.
Then, the Hilbert series of the graded algebra $\k[\Rep_\dd{A}]^{G_\dd}$
 is given by 
$$
h({\k[\Rep_\dd{A}]^{G_\dd}})=\Int(V,L,\dd)/\lambda(L\lo).
$$
\end{proposition}

\begin{proof}
Associated with the zero fiber of the map $p_\dd$, see
\eqref{flat}, is  the standard Koszul complex with differential $d_K$
of degree $-1$:
\begin{eqnarray}\label{BK}
&\BK_\idot(\Rep_\dd A)=\C[\VV_\dd ]\o\La^\hdot(\LS)^*
,\quad\text{such that}\\
&
\C[\Rep_\dd A]=\Coker[d_K: \BK_1(\Rep_\dd A)\to\BK_0(\Rep_\dd A)].\nonumber
\end{eqnarray}

For each $p\geq 0,$ the $\C^\times$-action on $\VV_\dd$ and on $\LS$
gives an additional  {\em weight grading}
on the vector space $\BK_p(\Rep_\dd A)$. 
The  Koszul differential respects the weight grading.
Writing $[r]$ to denote weight $r$ homogeneous component,
we get a bigraded direct sum decomposition
\beq{gradings}
\BK_\idot(\Rep_\dd A)=\bigoplus\nolimits_{p,r\geq0}  \BK_p(\Rep_\dd A)[r].
\eeq

Observe next that all the maps in \eqref{LS_exact}
are $G_\dd$-equivariant, provided one treats
$(L\lo)^*$ as a vector space with the trivial
$G_\dd$-action.
Therefore, dualizing  \eqref{LS_exact} yields a
 $G_\dd$-equivariant short exact sequence
\beq{psi}
0\too L\lo\stackrel{\psi^\top}\too
\LL_\dd^*\too (\LS)^*\too 0.
\eeq
Thus, one can rewrite the Koszul complex in the form
\beq{BK_equiv}
\BK_\idot(\Rep_\dd A)=\C[\VV_\dd] \o
\La^\hdot(\LL_\dd^*/L\lo).
\eeq

Assume now that the map
$p_\dd$ in \eqref{diag}
 is flat so that the scheme $\Rep_{\mathbf d}A$ is a complete intersection. 
Then, the Koszul complex provides a graded,
$G_\dd$-equivariant DG algebra resolution of
the coordinate ring $\k[\Rep_\dd{A}]$. 
It follows, since the group
$G_\dd$ is reductive, that the $G_\dd$-fixed part of  the Koszul complex provides a graded
 DG algebra resolution of $\k[\Rep_\dd{A}]^{G_\dd}$,
the subalgebra of $G_\dd$-invariants in  $\k[\Rep_\dd{A}]$.
Thus,
by the Euler-Poincar\'e principle, the Hilbert 
series of $\k[\Rep_\dd{A}]^{G_\dd}$ may be expressed in terms of
Hilbert series of the graded spaces $\BK_p(\Rep_\dd A)^{G_\dd}=
\bigoplus_{r\geq 0}\BK_p(\Rep_\dd A)^{G_\dd}[r]$ as follows
\begin{eqnarray}\label{BKh}
&h(\k[\Rep_\dd{A}]^{G_\dd})&=
\sum\nolimits_{p\geq 0} (-1)^p\cd
h\big(\BK_p(\Rep_\dd A)^{G_\dd}\big)\\
&&=
\sum\nolimits_{p\geq 0} (-1)^p\cd
h\left(\big(\C[\VV_\dd] \o
\La^p(\LL_\dd^*/L\lo)\big)^{G_\dd}\right).\nonumber
\end{eqnarray}

To compute the last expression, 
recall that for any two linear maps, $u:M\to M$ and
$v: N\to N$, of finite dimensional vector spaces, 
one has
$$
\sum_{p,m\ge 0}(-1)^p\cd t^m\cd \Tr(v\o
u|_{\sym^mN\o\La^pM})=\frac{\det(\bone-t\cd u)}{\det(\bone-t\cd v)}. 
$$

Now, for any $p\geq 0$, we have a direct sum decomposition with respect
to the weight grading
\beq{kkk}\C[\VV_\dd] \o
\La^p(\LL_\dd^*/L\lo)=
\bigoplus\nolimits_{r\geq 0}\big(\C[\VV_\dd] \o
\La^p(\LL_\dd^*/L\lo)\big)[r].
\eeq
Thus, using the above formula
and the notation $\bc_{ij}^r=\dim V_{ij}[r]-\dim L_{ij}[r]$, 
for any $g\in G_\dd,$ we find

\begin{eqnarray}\label{av}
&\sum\limits_{r,p\geq 0} (-1)^p  \cd t^r\cd
\Tr\left(g|_{\left(\C[\VV_\dd] \o
\La^p(\LL_\dd^*/L\lo)\right)[r]}\right)\qquad\qquad\qquad\\
&=\prod\limits_{r>0}\left({\Large
\frac{1}{(1-t^r)^{m_r}}\cdot\prod\limits_{i,j\in I}\frac
{\det(\bone-t^rg|_{\EE_{ij}})^{\dim L_{ij}[r]}}
{\det(\bone-t^rg|_{\EE_{ij}})^{\dim V_{ij}[r]}}}\right)=
\frac{1}{\lambda(L\lo)}\!\!&\!\!\prod\limits_{i,j\in I}
\frac{1}
{\det(\bone-t^rg|_{\EE_{ij}})^{\bc_{ij}^r}}.
\nonumber
\end{eqnarray}

Recall that, for any finite dimensional
$G_\dd$-module $M$, the dimension
of the $G_\dd$-fixed point subspace is given
by the integral
$\dim M^{G_\dd}=\int_{\UU_\dd} \Tr(g|_M)dg$.
We apply this to each direct summand 
of the graded $G_\dd$-module in the right-hand side 
of
 \eqref{kkk}.
We conclude that the  Hilbert  series of the subcomplex of
$G_\dd$-invariants from the last line in \eqref{BKh}  is obtained 
by averaging the expression in \eqref{av} 
over $\UU_\dd$, which gives the desired formula. 
\end{proof}

\subsec{Asymptotics of the matrix integral.}\label{asymptotic}
We will  use the following result about random unitary matrices  
(see e.g. \cite{DS}).  
 
\begin{proposition} \label{randommat} Let $m_s$, $n_s$, $s\ge 1$,  
be nonnegative integers, which are almost all zero.   
Then for large enough $N$, 
$$ 
\int_{U(N)}\;\prod_{s\ge 1}\big({\rm Tr}(u^s)\big)^{m_s}\big(\overline{{\rm Tr}(u^s)}\big)^{n_s}du= 
\prod_{s\ge 1}s^{m_s}\cd m_s!\cd \delta_{m_s,n_s} $$ 
\end{proposition} 
 
The statement of the proposition means that
the functions $z_s(u)={\rm Tr}(u^s)/\sqrt{s}$ behave, asymptotically as $N\to \infty$,   as  
independent complex-valued random variables distributed 
with Gaussian density $p(z)=\pi^{-1}e^{-z\overline z}$.  It follows, in
particular, that one has

\begin{corollary}\label{mat_cor} For any finite collection of polynomial functions 
$f_s\in\C[x_1,x_2], \,s=1,2,\ldots,m,$ there exists $N\gg0$ such that
for all $d\geq N$
one has
\beq{int1}\int_{U(d)}\;\left(\prod_{s\ge 1}^mf_s\big(z_s(u),\overline{z_s(u)}\big)\right)du=
\prod_{s\ge 1}^m\left(\int_{U(d)}\;f_s\big(z_s(u),\overline{z_s(u)}\big)du\right).
\eeq

Moreover,  we have
\beq{int2}
\int_{U(d)}\;f_s\big(z_s(u),\overline{z_s(u)}\big)du=
\frac{1}{\pi}\int_{\C}\;f_s(z,\bar z)e^{-|z|^2}dzd\bar z,
\quad\forall s=1,2,\ldots.
\eeq
\end{corollary}

The statement of Corollary \ref{mat_cor} can be generalized further
to  the case of integrals over $\UU_\dd=\prod_{i\in I}U(d_i)$,
a product of unitary groups. With such a generalization at hand, we are now ready to prove

\begin{corollary}\label{assymcor} For each integer $m$,
the sequence $\Int_m(V,L,\dd)$,
of $m$-th coefficients in the expansion \eqref{Int},
stabilizes as $\dd\to\infty$. Write $
\underset{\dd\to\infty}\lim\Int_m(V,L,\dd)$ for this stabilized  coefficient
and define a formal power series
$$\underset{^{\dd\to\infty}}\lim\Int(V,L,\dd):=\sum\nolimits_{m\geq 0}\, t^m\cd
\underset{^{\dd\to\infty}}\lim\Int_m(V,L,\dd).$$

Then, we have
$\dis\underset{^{\dd\to\infty}}\lim\Int(V,L,\dd)=\ze(V,L).$
\end{corollary}

\begin{proof} We have the standard identity
$$\det(\Id-u)=\exp{\big(\Tr\log(\Id-u)\big)}=
\exp\big(\sum\nolimits_{s\geq 1}\frac{1}{s}\Tr(u^s)\big)=
\prod_{s\geq 1}\exp\left(\frac{1}{s}\Tr(u^s)\right).$$

We apply this identity to the function under the integral sign in
formula \eqref{Int} and get
\beq{inte}
\prod_{r>0}\prod_{i,j\in I}\;
\frac{1}{\det(\bone_{\EE_{ij}}-t^r\cdot g_i^\vee\o
g_j)^{\bc_{ij}^r}}=
\prod_{s\ge 1}\exp
\biggl(\sum_{i,j\in I, r\in{\mathbb N}}
\bc_{ij}^r\cd\Tr(g^s|_{\EE_{ij}})\cd\frac{t^{rs}}{s}\biggr).
\eeq

Given $g=\{g_i\}_{i\in I}\in \UU_\dd,$ 
let $z_{is}(g)$ denote the function $g\mto z_{is}(g)=
{\rm Tr}(g_i^s|_{\EE_{ii}})/\sqrt{s}$.
According to formla
\eqref{Int}, using \eqref{inte}, we find
\beq{inte2}\Int(V,L,\dd;t)=\int_{U_\dd}\;
\left(\prod_{s\ge 1}\exp\big({\sum\nolimits_{i,j\in I,r\in{\mathbb N}}\bc_{ij}^r\cd t^{rs}\cd z_{is}(g) 
\overline{z_{js}(g)}}\big)\right).
\eeq

Now, using  an analogue of \eqref{int1} for the group $U_\dd$,
in the limit $d_i\to \infty$, $i\in I$, from \eqref{inte2} we obtain
$$ 
\lim_{\mathbf d\to \infty}
\Int(V,L,\dd;t)=\prod_{s\ge 1}\left(\int_{U_\dd}\;
\exp\big({\sum\nolimits_{i,j\in I,r\in{\mathbb N}}\bc_{ij}^r\cd t^{rs}\cd z_{is}(g) 
\overline{z_{js}(g)}}\big)\right)=
\prod\nolimits_{s\geq 1}
f_s( t), 
$$ 
where the functions $f_s$ are given by an  analogue of formula
\eqref{int2},
that is,
\begin{eqnarray*}
&f_s(t)&=\frac{1}{\pi^{|I|}}\int_{\mathbb C^I}\exp({\sum_{i,j\in
I,r\in{\mathbb N}}
\bc_{ij}^r\cd t^{rs}\cd  z_i
\overline{z_j}})\cd e^{-\sum_{k\in I} |z_k|^2}\\
&&=
\frac{1}{\pi^{|I|}}\int_{\mathbb C^I}\exp\big(-{\sum_{i,j\in
I,r\in{\mathbb N}}
(\delta_{ij}-\bc_{ij}^r\cd t^{rs})\cd  z_i
\overline{z_j}}\big)=
\frac{1}{\det(\bone-\sum_{r>0}\bc_{ij}^r\cd t^{rs})}.
\end{eqnarray*}

Finally, the last expression on the right equals
$1/\det \big(\bone-h(V; t^s)+h(L; t^s)\big),
$ by  definition of the coefficients  $\bc_{ij}^r:=\dim V_{ij}[r]-\dim L_{ij}[r]$.
\end{proof}

\begin{definition}\label{arci} 
A $(V,L)$-datum (or the algebra $A=T_RV/(L)$) is called an
{\em asymptotic RCI}
if, for any pair of positive integers $r,N,$ there exists a
dimension vector $\mathbf d\succ N$ such
that one has:
$$H_p\left(\bigoplus\nolimits_{s\leq r}\, \BK_\idot(\Rep_\dd
A)[s],\,d_K\right)=0,\quad\forall p>0,
$$
i.e., such that  all nonzero homology groups of 
the weight degree $\le r$ part  of the Koszul complex
 \eqref{gradings} vanish.
\end{definition}

Clearly, any RCI algebra is an  asymptotic RCI; we will see 
in \S\ref{Ex} below that
the class of RCI algebras is rather restrictive while
the class of   asymptotic RCI is much less restrictive.

Let $A$ be an  asymptotic RCI, and fix a positive integer $r.$
Combining together Proposition \ref{intform} and Corollary \ref{assymcor},
we deduce that there exists
a sequence of dimension vectors $\dd_1,\dd_2,\ldots,$
such that $\dd_m\to\infty$ and such that
the  integer valued function
$m\mto\dim(\k[\Rep_{\dd_m}{A}]^{G_\dd}[r])$
 stabilizes as $m\to\infty$ (we will see later, in 
Proposition \ref{surjinj},
that this stabilization is a general phenomenon that
holds for all large enough dimensions $\dd$ and, moreover,
has nothing to do with
the RCI property).
Write $\underset{^{\dd\to\infty}}\lim \dim(\k[\Rep_\dd{A}]^{G_\dd}[r])$
for  this stable value and define a power series
$$\underset{^{\dd\to\infty}}\lim h(\k[\Rep_\dd{A}]^{G_\dd}):=
\sum\nolimits_{r\geq 0}\, t^r\cdot
\underset{^{\dd\to\infty}}\lim \dim(\k[\Rep_\dd{A}]^{G_\dd}[r]).
$$
 
Combining together Propositions \ref{intform} and \ref{randommat},
we obtain the following result

\begin{proposition}\label{hrep} For any  asymptotic  RCI algebra
$A=T_RV/(L)$ we have, cf. \eqref{AA}$\mathsf{(ii)}$
$$
\underset{^{\dd\to\infty}}\lim h(\k[\Rep_\dd{A}]^{G_\dd})
=\frac{1}{\lambda(L\lo) }\cdot\underset{^{\dd\to\infty}}\lim \Int(V,L,\dd)=
\ze(V,L)/\lambda(L\lo).
$$
\end{proposition}

\section{{ {Noncommutative complete intersections}}}\label{sec2}
\subsec{NCCI algebras.}\label{NCCI}
Below, we consider
 nonnegatively graded algebras $A=\bigoplus_{r\geq 0}A[r]$
such that $A[0]=R$, to be referred as $R$-algebras.
We write $A_+=\bigoplus_{r>0}A[r]$,
so $A/A_+=R$.

For any $R$-bimodule $V$ and any two-sided ideal
$J\sset T^+_RV,$ the quotient $J/J^2$ has an obvious 
structure of  $T_RV/J$-bimodule.

Let $V$ be an $\N$-graded finite dimensional $R$-bimodule
and $J\sset T^+_RV$ a two-sided graded ideal.
A graded  algebra of the form $A=T_RV/J$  is called a 
{\em noncommutative complete intersection}
(NCCI) algebra provided the equivalent conditions (i)-(iii)
of the theorem below hold for $A$.

\begin{theorem}\label{ncci} Let $V$  be a  finite dimensional $\mathbb N$-graded
$R$-bimodule and $J\sset T_RV$ a two-sided finitely generated
graded ideal.
The following properties of the graded algebra
$A=T_RV/J$ are equivalent:

\vi $J/J^2$ is  projective as a graded $A$-bimodule;

\vii   \parbox[t]{135mm}{The algebra $A$ has
Hochschild dimension $\leq 2,$ i.e.,
one has $\Ext_{\bimod{A}}^i(A, M)=0,$ for any $A$-bimodule $M$ and any
$i\geq 3;$}

\viii  \parbox[t]{135mm}{The algebra $A$ has global  dimension $\leq 2,$ i.e.,
one has $\Ext_{\Lmod{A}}^i(A/A_+,M)=0,$ for any left
 $A$-module $M$ and any $i\geq 3.$}
\end{theorem}

This result is known; its proof will be recalled in \S\ref{Anick_proof}
below.

\subsec{Anick's resolution and the noncommutative Koszul complex.}
\label{kos}
Fix a $(V,L)$-datum and let $A=T_RV/(L).$

Following \cite{AH}, \cite{An}, \cite{GS}, one introduces a noncommutative
Koszul complex, $\K_\idot{A},$ (referred to as the
{\em Shafarevich complex} in \cite{Go}, \cite{Pi}), as follows.
Set $\K_\idot{A}=T_R(V\oplus L)$, and view this tensor algebra as
a graded algebra with respect to 
{\em homological grading}, cf. \S\ref{super},
such that the vector space $V\sset T^1_R(V\oplus L)$
is placed in  homological degree 0 and  the vector space $L\sset T^1_R(V\oplus L)$
is placed in  homological degree 1. Thus, we have
$\K_0{A}=T_RV.$ Define the map
\beq{d}
V\oplus L\too T_R(V\oplus L),\quad v \oplus \ell\mto j(\ell)\in T_RV\sset 
T_R(V\oplus L)
\eeq
Let
$d_K: \K_\idot{A}\to \K_{\hdot-1}{A}$  be the unique
super-derivation of the graded algebra $\K_\idot{A}$ 
that is given on the generators from $V\oplus L\sset\K_\idot{A}$
by formula \eqref{d}. It is clear that we have $d_K^2=0$.
Thus, $(\K_\idot{A}, d_K)$ is a DG algebra.
Write $H_\idot(\K_\idot{A}, d_K)$ for the corresponding homology.
Clearly we have $H_0(\K_\idot{A}, d_K)=A$. Therefore,
the total  homology space $H_\bullet(\K_\idot{A}, d_K)$ 
acquires the natural  structure of a graded $A$-algebra.

Next, we recall the definition of 
Anick's resolution \cite{An2}. This is a complex  of graded $A$-bimodules
of
the form
\beq{C}
C_\idot{A}:\en
A\otimes_R L\otimes_R A\stackrel{\pa}\too A\otimes_R V\otimes_R A
\stackrel{f}\too A\otimes_R
A\stackrel{\bm_A}\too A\to 0.
\eeq

Here, $\bm_A$ is the multiplication map,
and the maps $\pa$ and $f$ are given by the formulas
$$
\pa(a_1\otimes \ell\otimes a_2)=
a_1\cd D(\ell)\cd a_2,
\quad\text{resp.}\quad
f(a_1\otimes v\otimes a_2)=a_1v\otimes a_2-a_1\otimes va_2.
$$
In the first formula above, we have used the map 
$$D: T_R^+V\to A\otimes_RV\otimes_R A,
\en
v_1\otimes\ldots\otimes v_n\mto
\sum_{p=1}^n (\overline{v_1\otimes\ldots\otimes v_{p-1}})\otimes 
v_p\otimes (\overline{v_{p+1}\otimes\ldots\otimes v_n}),
$$
where bar stands for the image  in $A$ of an element of
$T_RV$.

\begin{definition}\label{min} A $(V,L)$-datum (or a graded $R$-bimodule $L\sset
T_R^+V$)
is said to be
{\em minimal} provided one has  
$$(T_RV\cdot L\cdot T_R^+V+ T_R^+V\cdot
L\cdot T_RV)\cap L=0.
$$
\end{definition}

It is clear that, for any graded  $R$-bimodule  $L\sset T^+_RV$, one can always find
a minimal sub-bimodule $\widetilde{L}\sset L$ such that
$(L)=(\widetilde{L}).$ Although the choice of such a minimal
 sub-bimodule $\widetilde{L}$ is not unique, any
two   minimal sub-bimodules $\widetilde{L}$ have the same Hilbert
series.

For any  $(V,L)$-datum,
the composite $L\into J\onto J/J^2$ extends, by $A$-linearity,
to an $A$-bimodule map $\pi: A\o_R L\o_R A\to J/J^2$.

The proof of the following known (see e.g. \cite{An})
result will be given in \S\ref{Anick_proof}.

\begin{theorem}\label{Anick} The complex $C_\idot{A}$ is exact in all the terms
except possibly the first one. Further,
 the following conditions are equivalent

$\mathsf{(1)}\quad$  The bimodule $L$ is minimal and $A$ is an NCCI.

$\mathsf{(2)}\quad$  The  Hilbert  series of $A$ equals 
$\dis
h(A)=\big(\bone-h(V)+h(L)\big)\inv;
$

$\mathsf{(3)}\quad$  The above defined natural map $\pi:
A\o_RL\o_R A \to J/J^2$ is an isomorphism;

$\mathsf{(4)}\quad$  \parbox[t]{130mm}{The complex $C_\idot{A}$ is an 
$A$-bimodule resolution of $A$, i.e.
the first map $\pa$ in \eqref{C} is injective;}

$\mathsf{(5)}\quad$ \parbox[t]{130mm}{The Koszul complex $\K_\idot A$ is a
DG algebra resolution of $A$, i.e., we have $H_k(\K_\idot{A},d)$
$=0,$ for any  $k >0$.}
\end{theorem}
We remark that if condition (5) holds then $K_\idot{A}$ is, in effect,
a {\em minimal} free associative DG algebra resolution of $A$ in the sense of
minimal models, cf. \cite{AH}.

\subsec{Noncommutative differential forms.}\label{NCOm}
Associated with an $R$-algebra $A$, is an $A$-bimodule
 $\Om^1A$ of {\em noncommutative} relative
1-forms on $A$ defined as the kernel of the multiplication map
$\bm_A: A\o_RA\to A.$ Thus,   one has a
 short exact sequence of $A$-bimodules
\beq{fund}
0\too
\Om^1A\stackrel{\imath_A}\too A\o_RA\stackrel{\bm_A}\too A\too 0.
\eeq

Next, let $F$ be an $R$-algebra,  $J\sset F$ a two-sided ideal, and
let  $A:=F/J$. We put $\Om^1(F|A):=A\o_F\Om^1F\o_F A$.
This is an $A$-bimodule  equipped with an
$A$-bimodule map $\tau: \Om^1(F|A)\to \Om^1A$,
induced by the projection $F\onto A$.
There is a canonical short exact sequence
of $A$-bimodules, cf. \cite[Corollary 2.11]{CQ},
\beq{CQ}
0\too J/J^2 \stackrel{d}\too
\Om^1(F|A)\stackrel{\tau}\too \Om^1A\too 0.
\eeq
Here, the  map 
$J/J^2\to
A\o_F\Om^1F\o_FA=\Om^1(F|A)$ is induced by  restriction to $J$ of
the de Rham differential $d: F\to \Om^1F$.
This  exact sequence  may be thought
of as a noncommutative analogue of
the {\em conormal exact sequence} of a subvariety.

We now fix a $(V,L)$-datum and set
 $F=T_RV,\,J=(L)$, and $A=T_RV/(L)$.
Then we have  natural isomorphisms
\beq{Om}
\Om^1F\cong F\otimes_RV\otimes_R F,
\quad\text{and}\quad\Om^1(F|A)= A\otimes_RV\otimes_R A.
\eeq

Splicing \eqref{CQ} and \eqref{fund} together, we obtain
the exact sequence in
the bottom row of the following commutative
diagram of $A$-bimodules
\beq{diag}
\xymatrix{
\Ker(\pa)\ar@{^{(}->}[r]& A\otimes_R L\otimes_R A\en\ar[r]^<>(0.5){\partial}
\ar@{->>}[d]_<>(0.5){\pi}&
A\otimes_R V\otimes_R A\ar[r]^<>(0.5){f}\ar@{=}[d]^<>(0.5){\eqref{Om}}&
A\otimes_RA\ar[r]^<>(0.5){\bm_A}\ar@{=}[d]&A\to 0\\
0\ar[r]&
J/J^2\en\ar[r]^<>(0.5){d}&\Om^1(F|A)\ar[r]^<>(0.5){\imath_A\ccirc\tau}
&A\o_R A\ar[r]^<>(0.5){\bm_A}&A\to0.
}
\eeq

\subsec{Proof of Theorem \ref{ncci} and   Theorem \ref{Anick}.}\label{Anick_proof}
 The first claim of Theorem \ref{Anick}
amounts to the exactness  of the top row of diagram
\eqref{diag}. The latter follows from  commutativity of the
 diagram since the bottom row of the diagram is
 exact.
Also, from the commutative diagram
 we deduce $\Ker\pa\cong \Ker\pi.$ This
yields  the equivalence (3) $\Leftrightarrow$~(4).

Next, we apply the Euler-Poincar\'e 
principle to the complex
$C_\idot A$ (in which the differential has degree $0$), 
and find 
$$
h(A) \cd\Big(\bone-h(A)\cd \big(\bone-h(V)+h(L)\big)\Big)+h({{\rm Ker}\pa})=0. 
$$
Thus $h({\Ker \pa})=0$ if and only if  $h(A)=(\bone-h(V)+h(L))^{-1}$. 
This yields the equivalence (2) $\Longleftrightarrow$ (3).
Also, it is clear that
(2) implies property (i) of Theorem \ref{ncci}.

Now, if  $J/J^2$  is a  projective $A$-bimodule,
then the bottom row of diagram \eqref{diag} provides a
length 3 projective $A$-bimodule resolution of $A$. It follows that
$A$ has Hochschild dimension $\leq 2$.
Hence, in Theorem \ref{ncci}, we have (i) $\Rightarrow$ (ii).
Similarly,  we have (i) $\Rightarrow$ (iii).
Conversely, from the exact sequence in the second row of diagram
\eqref{diag},
 for any $A$-bimodule $M$, we deduce
$\Ext_{\bimod{A}}^1(J/J^2,M)=\Ext_{\bimod{A}}^3(A, M).$
The latter Ext-group vanishes if 
$A$ has Hochschild dimension $\leq 2$. This yields the 
implication  (ii) $\Rightarrow$ (i)
 in Theorem \ref{ncci}.

Any {\em graded} finitely generated
projective  $A$-bimodule has the form $A\o_R\widetilde{L}\o_RA$,
for some  finitely generated graded $R$-bimodule $\widetilde{L}$.
In particular, if  $J/J^2$ is projective, we deduce
$$(A/A_+)\o_A(J/J^2)\o_A(A/A_+)=(A/A_+)\o_A(A\o_R\widetilde{L}\o_RA)\o_A(A/A_+)=
\widetilde{L}.
$$
Therefore,
applying the functor $(A/A_+)\o_A(-)\o_A(A/A_+)$ to the vertical
map $\pi$ in diagram \eqref{diag}, we obtain a surjective
 $R$-bimodule map
$$L=(A/A_+)\o_A(A\o_R L\o_RA)\o_A(A/A_+)\,\onto\,
(A/A_+)\o_A(J/J^2)\o_A(A/A_+)=
\widetilde{L}.
$$

The resulting map $L\to \widetilde{L}$ may be identified with
the composite map 
$$
L\into J\onto
\frac{J}{F_+\cd J+J\cd F_+}=\frac{J}{J^2+F_+\cd J+J\cd
F_+}=(A/A_+)\o_A(J/J^2)\o_A(A/A_+).
$$
Observe that the bimodule $L$ is minimal if and only if the latter
map  is a bijection. Since $L$ is minimal by assumptions,
we deduce that the map $\pi$ induces  an  isomorphism
 $L\iso \widetilde{L}.$ It follows that
 the vertical
map $\pi$ in diagram \eqref{diag} is itself  an  isomorphism,
and we have proved the equivalence (1) $\Longleftrightarrow$ (3).

The equivalence (1) $\Longleftrightarrow$  (5), as well
as the implication  (i)  $\Rightarrow$ (iii) of Theorem \ref{ncci},
are due
to Anick, see  Theorem 2.6 and Corollary~2.12(b) in \cite{An};
cf. also  Remark 5.3.5 in \cite{Gi3} for an alternative approach.
\qed
\vskip 3pt

The implication (5) $\Rightarrow$  (2) has been observed first  by 
Golod-Shafarevich \cite{GS}.

\subsec{A super-version}\label{super}
Given a $\Z/2\Z$-graded {\em super}-vector space
 $M=M_\even\oplus M_\odd$,
for the corresponding {\em super}-symmetric
(supercommutative) algebra we introduce the notation
\beq{ssym}
\ssym M:=(\sym M_\even)\;\otimes\;(\La M_\odd).
\eeq
If each of the spaces
$M_\even, M_\odd$ is equipped with an additional
grading, we say that $M$ is a graded
super-vector space. In that case, we view $\ssym M$ as a graded
super-algebra equipped with the {\em total} grading.

Below,
we will frequently encounter  {\em bi}-graded vector spaces
with two types of gradings of the form
$M=\bigoplus_{p\in\Z,r\geq 0} M_p[r]$.
The  $p$-grading, written
as $M=\bigoplus_{p\in\Z} M_p$
where  $M_p:=\bigoplus_{r} M_p[r]$, 
 will be called
{\em homological} grading.
The   $r$-grading, written
as $M=\bigoplus_{r\geq 0} M[r]$ where    $M[r]=\bigoplus_{p} M_p[r],$
 will be called {\em weight} grading.
In such a case, we put
$M_\even=\bigoplus_{p\;\text{even}} M_p$ and
$M_\odd=\bigoplus_{p\;\text{odd}} M_p.$
The bigrading on  $M$ induces corresponding homological
and weight  gradings $\ssym M=$
$\bigoplus_{p\in\Z,r\geq 0}(\ssym M)_p[r]$. Thus, for any $k\geq 0$, we have
$$m_1\& \ldots \& m_k\in (\ssym M)_{p_1+\ldots+p_k}[r_1+\ldots+r_k],
\quad\forall
m_1\in M_{p_1}[r_1],\ldots, m_k\in M_{p_k}[r_k].
$$

Hilbert series will be always
taken with respect to the weight grading
(with parity being determined by the homological grading).
Specifically, we define  the  Hilbert  series 
of a bi-graded (super)-vector space $M=\bigoplus_{p\in\Z,r\geq 0} M_p[r]$ to be
\beq{super_h}
h(M; t)=h(M_\even;t)-
h(M_\odd;t)=\sum_{p,r\geq 0} (-1)^p\cdot t^r\cdot \dim  M_p[r].
\eeq

Let $F=\bigoplus_{r\geq 0} F[r]$ be a graded
super-algebra with $R=F[0]\sset F_\even$.
Set $F_+=\bigoplus_{r> 0} F[r]$.
We define $[F,F]_{\text{super}}$,
the super-commutator space,  to be the linear span 
of  elements of the form $ab-(-1)^{|a|\cdot |b|}ba$, where 
$a,b$ run over even and odd  elements of $F$,
and we write $|x|=0$ if $x\in F_\even$,
resp. $|x|=1$ if $x\in F_\odd$.
 Observe that
$[F,F]_{\text{super}}\sset F_+$, since $R$ is a commutative
algebra.

We extend the notation
from \eqref{OA} to the super-algebra setting and put, cf. \eqref{ssym},
\beq{OAsuper}
\O(F):=\ssym\Big(F_+\big/[F,F]_{\text{super}}\Big).
\eeq
The weight, resp. homological, grading on $F$ descends to the super-commutator
quotient and makes $\O(F)$ a super-commutative graded, or bigraded, super-algebra,
$\O(F)=\bigoplus_{r\geq 0}\O(F)[r],$
with the corresponding homological, resp. weight,
gradings.

An ordinary algebra $A$ may  be treated as a  super-algebra with
 $A_\odd=0$.
 In this case,  formula  \eqref{OAsuper} reduces
to formula \eqref{OA}, i.e., we have
$\O(A)=\sym\big(A_+\big/[A,A]\big)$.

The results of previous sections can be generalized to the case
where each of the $R$-bimodules in the $(V,L)$-datum
 is a super-vector space.
Thus, in such a case, we have
$V=V_\even\oplus V_\odd,\, L=L_\even\oplus L_\odd$
and, moreover, the imbedding $j: L\into T^+_RV$
respects parity, i.e., is a morphism of super-vector spaces.
Then, we define the  super-space $L\lo:=j\inv([T_RV,T_RV]\super)$, and set 
$m_r=\dim (L\lo)_\even[r]-\dim (L\lo)_\odd[r]$.

The theory of representation schemes in the supercase is parallel
to the even case, except that these schemes are now not ordinary affine
schemes but rather affine superschemes (so the rings $\k[\Rep_\dd A]$
are supercommutative). The notion of NCCI algebra 
is defined analogously to the even case. The notion of RCI algebra
cannot be based on a dimension equality
like \eqref{dim}, since
dimension of singular super-schemes is ill-defined, in general.
Instead, one defines  RCI algebras by requiring that the
map $p_\dd$ in diagram \eqref{diag} be flat which, in the even case, is equivalent
 to \eqref{dim} by Lemma \ref {meaning}.

\subsec{Karoubi-de Rham complex and cyclic homology.}\label{HC} 
We write $\HC_\idot(A)$ for {\em  reduced} cyclic homology of an $R$-algebra $A$.
We will also use  reduced Hochschild
 homology of $A$, to be denoted 
 $\overline{HH}_\idot(A)$. 
For an augmented $R$-algebra
$A$,  one has $\overline{HH}_0(A)=\HC_0(A)=A\big/(R+[A,A]),$
and $HH_p(A)=\overline{HH}_p(A)$,
for all  $p>0$.

Let $\Om^\hdot A:=T_A(\Om^1A)$ be the DGA
of noncommutative differential forms, see \cite{CQ}.
The {\em reduced Karoubi-de Rham complex} of $A$ is defined as 
the super-commutator quotient 
$$\DR  A:=\Om^\hdot A\Big/\big(R+[\Om^\hdot A,\,\Om^\hdot A]_{\text{super}}\big).
$$
The natural  differential $\Om^\hdot A\to\Om^{\hdot+1}A$
descends to a  de Rham differential on $\DR  A$.

Further, according to \cite{Gi2}, there is a canonical
isomorphism 
$$
\overline{HH}_\idot(A)\cong \Ker[\imath:\DR^\hdot A\to\Om^{\hdot-1} A],
$$
where $\imath$ is a certain  canonical map that  anti-commutes with the
 differentials. 
The de Rham differential gets transported,
via the above isomorphism, to a natural differential
$d_H: \overline{HH}_\idot(A)\to \overline{HH}_{\idot+1}(A)$ which is
essentially (up to a nonzero
constant factor) induced by   Connes' differential $B$.

\begin{lemma}\label{hclem} Let $A$ be an $R$-algebra with vanishing
reduced de Rham homology, i.e. such that $H_\idot(\DR  A,d_{DR})=0$.
Then,  the following complex is an exact sequence
\beq{HHH}
0\too R\too {HH}_0(A)\stackrel{d_H}\too {HH}_1(A)\stackrel{d_H}\too {HH}_2(A)\stackrel{d_H}\too
{HH}_3(A)\too\ldots.
\eeq

Furthermore, there are isomorphisms
$$
\HC_j(A)\cong \im[B:HH_j(A)\to HH_{j+1}(A)],\en\;\forall j\geq 0.
$$
If, in addition, Hochschild homology groups vanish
in all degrees $j>m$, for some integer $m>1$, then we deduce
$\HC_{m-1}(A)\cong {HH}_m(A),$ and $\HC_j(A)=0$  for all $j\geq m$.
\end{lemma}

\begin{proof} There is a standard exact sequence of reduced homology
groups, 
see \cite[Thm. ~2.6.7]{Lo}:
$$0\to H_j(\DR {A}) \to \overline{HC}_j(A)
\stackrel{B}\too \overline{{HH}}_{j+1}(A).
$$
It follows that Connes' differential $B$ is injective,
 for any algebra with vanishing reduced de Rham homology.

Next, we use  Connes' long exact sequence, see \cite[\S2.2.13]{Lo},
\beq{connes}
\ldots\to\overline{HC}_{j+1}(A)\to
\overline{HC}_{j-1}(A)\stackrel{B}\too\overline{{HH}}_j(A)
\to\overline{HC}_j(A)\to\ldots.
\eeq
From the injectivity of $B$, we deduce that this  long exact sequence
breaks up into short exact sequences of the form
\beq{connes2}
0\too\overline{HC}_{j-1}(A)\stackrel{B}\too\overline{{HH}}_j(A)
\too\overline{HC}_j(A)\too0.
\eeq
Splicing all these  short exact sequences together yields
a long exact sequence as in \eqref{HHH}.
It is easy to verify that the maps in the two  sequences
are proportional, up to nonzero constant 
factors. Hence, the sequence in \eqref{HHH} is exact.

To complete the proof, assume that ${HH}_j(A)=0$ for all $j>m>1$.
Then, $\overline{HH}_j(A)=0$,
and we deduce that $\overline{HC}_k(A)=0$ for all $k\geq m$, due to
injectivity of $B$. Hence,
 $\overline{HC}_{m+1}(A)=\overline{HC}_m(A)=0.$
Now, using the  short exact sequence \eqref{connes2} for
$j=m,$ we get an isomorphism
$\HC_{m-1}(A)\iso \overline{HH}_m(A)=HH_m(A),$ and the Lemma follows.
\end{proof}

\subsec{Cyclic Koszul complex.}\label{CCI}
From now on, we fix a {\bf{minimal}}  $(V,L)$-datum. 
 We
define the {\em cyclic
Koszul complex} of the algebra $A=T_RV/(L)$
to be the {\em super}-commutator quotient 
$$
\KK_\idot{A}:=\K_\idot{A}\big/[\K_\idot{A},\K_\idot{A}]_{\text{super}},
$$

This is a bigraded vector space, $\KK_\idot{A}=\bigoplus_{p,r\geq
0}(\KK_pA)[r]$, and
the Koszul differential on $K_\idot{A}$  descends to 
a well-defined differential
$d_K: (\KK_\idot{A})[r]\to (\KK_{\hdot-1}A)[r]$.
Furthermore, we clearly have a graded
space isomorphism $H_0(\KK_\idot{A}, d_K)=A/[A,A].$

\begin{proposition}\label{HC1} Let $A=T_RV/(L)$ be an NCCI algebra.
Then, we have
\beq{h12}
 H_1\big(\KK_\idot{A},d\big)\cong {HH}_2(A),\quad
\text{and}\quad H_j\big(\KK_\idot{A},d\big)=0,\quad\forall j>1.
\eeq

Furthermore, for Hilbert series one has the formulas
\beq{ncci_h}
h\big(\O(K_\idot{A})\big)=\ze(V,L),
\quad\text{and}\quad h\big(\O(A)\big)=\ze(V,L)\cd h\big(\sym
{HH}_2(A)\big).
\eeq
\end{proposition}

\begin{proof} If $A$ is an NCCI algebra, then
 the DG algebra $\K_\idot{A}$ provides
a free resolution of $A$,
by Theorem \ref{Anick}(5).  Thus, according to the definition of cyclic
homology, the  complex $\big((\KK_\idot{A})/R,d_K\big)$ computes
reduced cyclic homology of the algebra $A$.

Now, we have that $A=\bigoplus_{r\geq 0} A[r]$ is a nonnegatively graded algebra.
It follows by Poincar\'e lemma
that the algebras $A$ and $A[0]$ have the same  Karoubi-de Rham
homology. Clearly,
 all  Karoubi-de Rham homology groups of the algebra $A[0]=R$
 vanish. Moreover, since $A$ is a NCCI algebra,
we have  ${HH}_j(A)=0$, for all $j>2.$ Thus, 
by Lemma \ref{hclem}, we conclude that
 $\HC_j(A)=0$, for all $j>1$ and, moreover,
 there is an isomorphism 
$\HC_1(A)\cong  {HH}_2(A)$.

Now,   view the noncommutative Koszul complex $K_\idot{A}$
as a free graded super-algebra,  $T_R(V\oplus L),$
such that the vector space
$V$ is even and 
the vector space
$L$ is odd. For such
 a free super-algebra,
the corresponding Hilbert series of 
$\O\big(T_R(V\oplus L)\big)$ can be
computed purely combinatorially; in accordance with formula
\eqref{AA} for a free (super)-algebra, one finds
$h\big(\O\big(T_R(V\oplus L)\big)\big)=\ze(V,L),$ cf. \cite{St}.
This yields the first formula in \eqref{ncci_h}.

Alternatively, the same formula may be deduced from 
general results of \S\ref{ev_grad} below, applied to the
super-algebra $T_R(V\oplus L)$. Specifically,
since any  free super-algebra is clearly RCI, 
the
above formula for the Hilbert series of $\O\big(T_R(V\oplus L)\big)$
is a special case of equation \eqref{oa}.
The proof of that equation is independent of
the intervening material.

Next, given a $\Z_+$-graded super-vector space
$M$, write $\lll M\rrr :=\lll M_\even\rrr -\lll M_\odd\rrr$
for the class of $M$
in the Grothendieck group of $\Z_+$-graded super-vector spaces.
Applying the Euler-Poincar\'e
principle to the complex $\KK_\idot{A}$ and using \eqref{h12}, we deduce
$$\lll\KK_\idot{A}\rrr= \lll H_0\big(\KK_\idot{A}\big)\rrr-\lll H_1\big(\KK_\idot{A}\big)\rrr
=\lll A/[A,A]\rrr-\lll {HH}_2(A)\rrr.
$$
Hence, for the corresponding symmetric algebras, we find
$$
h\big(\O(A)\big)=h\left(\sym\frac{A_+}{[A,A]}\right)=
h\big(\sym(\KK_\idot{A})_+\big)\cd h\big(\sym {HH}_2(A)\big).
$$

The last equation, combined with the first formula in \eqref{ncci_h}, 
implies the second formula in \eqref{ncci_h} 
since we have  $\O(K_\idot{A})=\sym\big((\KK_\idot{A})_+\big),$
by definition.
\end{proof}
\begin{remark}
 One may use the standard identity $h(\sym M)\cdot h(\La M)=1,$
where $\La M$ is viewed as a {\em super}-algebra,
to rewrite the second formula in \eqref{ncci_h}
in the following form, more compatible with \eqref{AA}(ii):
$$h\big(\O(K_\idot{A})\big)=\ze(V,L)\big/h\big(\La {HH}_2(A)\big).
$$
\end{remark}
\smallskip

The assignment $\ell\mto 0\oplus\ell$ gives an imbedding
$L \into V\oplus L= T^1_R(V\oplus L)=\K_1{A}.$
Therefore, the vector space $L\lo\sset L$, see \eqref {lo},
 may be viewed as a subspace 
of $\K_1{A}$.
By definition of $L\lo$ and of the differential  on
$\K_\idot{A}$, we have $d_K(L\lo)\sset [T_RV,T_RV]\sset 
 [\K_\idot{A},\K_\idot{A}]_{\text{super}}$.
Therefore, the image of $L\lo$ in $\KK_\idot{A}$ is annihilated by
the differential in the cyclic Koszul complex.
Thus, we obtain a canonical linear map
\beq{KK}
L\lo \to  H_1(\KK_\idot{A}, d_K).
\eeq

\begin{lemma}\label{easy} 
The map \eqref{KK} is injective. Furthermore, we have
$$H_j\big(\KK_\idot{A}/L\lo,\, d_K\big)=0,\en\forall j>0\quad
\Longleftrightarrow\quad
H_j(\KK_\idot{A},\, d_K)=
\begin{cases} L\lo&\operatorname{if}\en j=1\\
0&\operatorname{if}\en j>1.
\end{cases}
$$
\end{lemma} 
\begin{proof}
Both claims  follow
from the  fact that $L\lo\cap \im(d_K)=0,$ which is
 easily verified.
\end{proof}

The second important result of this paper,
to be proved in \S\ref{sec3}, reads

\begin{theorem}\label{converse} Fix  a  minimal $(V,L)$-datum and set
 $A=T_RV/(L)$.
 Then, the following
conditions \vi and \vii are equivalent:

\vi The algebra $A$ is NCCI  and the map \eqref{KK} is a bijection.

\vii The algebra $A$ is an asymptotic RCI, cf. Definition \ref{arci}.

Furthermore, the above conditions imply that 
 ${HH}_2(A)\cong L\lo$, and one has
$$h(\O(A))=\ze(V,L)/\lambda(L\lo),\quad\text{and}\quad h(\sym
{HH}_1(A))=\ze(V,L)/\lambda(L\lo)^2.
$$
\end{theorem}

The above formulas for Hilbert series are
identical to those in \eqref{AA}.
Since any RCI algebra is an  asymptotic RCI,
we see that
Theorem \ref{converse} implies Theorem \ref{maint}
(the long exact sequence in \eqref{exact}  follows
from that of Lemma \eqref{hclem}).
 
\section{The trace map}\label{sec3}
\subsection{Evaluation homomorphism.}\label{ev_sec}
Fix  an $R$ super-algebra $A$ and 
a  dimension vector $\dd$.
Thus, one has the super-scheme $\Rep_\dd A$.

Recall the notation $\EE_\dd =\End_\C(\bigoplus_{i\in I}\C^{d_i})$,
and write
$\EE_\dd \o\k[\Rep_\dd{A}]$, a 
 tensor product of two  associative algebras.
The group $G_\dd$ acts naturally on each tensor factor
and we let $\big(\EE_\dd \o\k[\Rep_\dd{A}]\big)^{G_\dd}$
denote the subalgebra of $G_\dd$-invariants with
respect to the diagonal action.

To each element $a\in A$, one associates
the  function $\wh{a}: \rda\to\EE_\dd ,$
$\rho\mapsto \wh{a}(\rho):=\rho(a)$. The assignment
$a\mapsto\wh{a}$ clearly gives an algebra
homomorphism, called  {\em evaluation map},
\beq{tr_om} \ev_\dd:\ A\too\big(\EE_\dd \o\k[\Rep_\dd{A}]\big)^{G_\dd},
\quad a\mto\wh{a}.
\eeq

Further, we have  the  (super)-trace map 
$$
\EE_\dd \o\k[\Rep_\dd{A}]
\xymatrix{\ar[rr]^<>(0.5){\Tr\o\Id}&&}
\C\o\k[\Rep_\dd{A}]=\k[\Rep_\dd{A}].$$
This map clearly  vanishes on
(super)-commutators. Therefore, the
composite of the homomorphism \eqref{tr_om} with the trace map above
 descends
to the (super)-commutator quotient space
of the algebra $A$. 
Thus, we obtain a well defined linear map
$A/[A,A]\super\to \k[\Rep_\dd]^{G_\dd}, \, a\mto \Tr(\wh{a}).$

Next, we assume that there is  an augmentation, $A\onto R$,
with augmentation ideal $A_+$.
The
 linear map $a\mto \Tr(\wh{a})$ may be uniquely  extended, by multiplicativity,
to a super-algebra  homomorphism 
\begin{eqnarray}\label{trace}
&\TR:\
\underset{}{\O(A)=\ssym\left({A_+\big/[A,A]\super}\right)}\too\k[\Rep_\dd]^{G_\dd},
\\
&a_1\&\ldots\& a_m\mto\Tr(\wh{a}_1)\cd\ldots\cd\Tr(\wh{a}_m),\quad
\forall
a_1\&\ldots\& a_m\in\ssym^m\left({A_+\big/[A,A]\super}\right).\nonumber
\end{eqnarray}

\begin{proposition}\label{surj} The map
$\TR$ is  surjective.
\end{proposition}
\begin{proof}
The claim follows from Weyl's fundamental theorem of invariant
theory; the ring of invariants of a collection of tensors 
is generated by various contractions of these tensors,
cf. \cite{LBP}.
\end{proof}

As has been noticed in \cite{Gi1}, the map $\TR$
 becomes `asymptotically bijective',
in a sense,
as $\dd\to\infty$. One way to make this heuristic idea precise
is to consider the case of algebras
equipped with an additional {\em weight} grading, as we are
going to do below.

\subsec{Trace map for graded algebras.}\label{ev_grad}
Fix a $(V,L)$-datum of graded $R$ super-bimodules
$V=V_\even\oplus V_\odd, \, L=L_\even\oplus L_\odd,$
and set $A=T_RV/(L)$. 

For any dimension vector $\dd$,
the  grading on $V$ gives rise
to a weight grading on the super-algebra
$\k[\Rep_\dd A]^{G_\dd}$, cf. \eqref{weight}.
%
%
%

The map $\TR$ in Proposition \ref{surj}
 clearly respects the gradings.

\begin{definition}\label{asbij} Assume that, for each dimension vector $\dd$,
we are given
a pair of graded vector spaces, $M_\dd=\bigoplus_{r\geq 0}M_\dd[r]$
and $N_\dd=\bigoplus_{r\geq 0}N_\dd[r]$, and a linear
map $f_\dd: M_\dd\to N_\dd$ that respects the gradings.

We say that the maps $f_\dd$ are
{\em  asymptotically bijective} as $\dd\to\infty$ provided,
for any  $r=1,2,\ldots,$ there exists
a positive integer $n(r)\gg 0$ such that
the  map 
$$
f_\dd:\ \bigoplus\nolimits_{s\leq r} M_\dd[s]
\too \bigoplus\nolimits_{s\leq r} N_\dd[s]
$$
 is a bijection for all $\dd\succ n(r).$
\end{definition}

Proposition 11.1.1 from  \cite{CBEG} yields the following result
\begin{proposition}[Stabilization]\label{surjinj} The maps
$\TR$ are 
 asymptotically bijective.\qed
\end{proposition}

Hence, 
for   the Hilbert series
$h\big(\k[\Rep_\dd{A}]^{G_\dd}\big),$ we deduce

\begin{corollary} \label{lim}
With the notation used in Corollary \ref{hrep}, we have 
$$
h\big(\O(A)\big)=\underset{^{\dd\to\infty}}\lim
h\big(\k[\Rep_\dd{A}]^{G_\dd}\big).\qquad\Box
$$
\end{corollary}

We combine the equation of Corollary \ref{lim} with the formula for
the asymptotics of the matrix integral studied in
\S\ref {asymptotic}. We conclude that,
for any  asymptotic RCI super-algebra $A=\bigoplus_{r\geq 0}A_p[r]$,
one has
\beq{oa}
h\big(\O(A)\big)\;
\xymatrix{
\ar@{=}[rr]^<>(0.5){\text{Coroll. } \ref{lim}}&&
\;\underset{^{\dd\to\infty}}\lim
h\big(\C[\Rep_\dd{A}]^{G_\dd}\big)\;\ar@{=}[rr]^<>(0.5){\text{Prop. } \ref{hrep}}&&
}\;\ze(V,L)/\lambda(L\lo).
\eeq

This proves the formula  (as well as its super-analogue)
for the Hilbert series of
$\O(A)$ stated in Theorem \ref{maint} and in Theorem \ref{converse}.

\begin{remark}
 Any free tensor super-algebra of the form $A=T_RV$ is
clearly an RCI algebra, associated with the $(V,L)$-datum such that
 $L=0$.
In this special case, the vector space
$A/[A,A]\super$ has a basis of cyclic words (with the
understanding that the cyclic word $w^{2m}$ is equal to zero,
for any word
$w$ of odd homological degree and any $m\geq 1$), and 
there is an alternative purely algebraic proof of 
formula \eqref{oa},
which does not use matrix integrals
  and holds over fields of any
characteristic, cf. \cite{St}.  
\end{remark}

\begin{example} Let $A=\C[x]$, where $x$  is odd; so the
 algebra
$A$ has no  relations.
In this case $[A,A]\super$ is spanned by
$x^2,x^4,\ldots$, so $A_+/[A,A]\super$ has basis $x,x^3,x^5,..$ and is purely
odd. Thus $\O(A)=\La(x,x^3,\ldots),$ and we easily find
$$
h({\O(A)}; t)=(1-t)(1-t^3)(1-t^5)\cd\ldots.
$$

On the other hand, the super-version of the formula
from Theorem \ref{maint}, equivalently, formula
\eqref{oa},  claims that 
$$
h({\O(A)}; t)=[(1+t)(1+t^2)(1+t^3)\cd\ldots]^{-1}. 
$$
Thus we deduce  the following well known classical identity
$$
(1-t)(1-t^3)(1-t^5)\cd\ldots=[(1+t)(1+t^2)(1+t^3)\cd\ldots]^{-1}.
$$
\end{example}

\subsec{Adding dummy variables.}\label{dummy_sec} The technique of 
`dummy variables' exploited below is not new, it has been successfully
applied in similar circumstances  by C. Procesi, and others,
cf. \cite{Pr}.

Given a finite dimensional $R$-bimodule $V$,
the vector space $V^*=\Hom_\C(V,\C)$
 comes equipped with a natural $R$-bimodule structure,
that is, with a bigrading
$V^*=\bigoplus_{i,j\in I} V^*_{ij}$
such that $V^*_{ij}=(V_{ji})^*.$

Now, fix a $(V,L)$-datum and put $V\z:=V\oplus V^*\oplus L^*,$
where $V\z$ is viewed as a  $R$-{\em  super}-bimodule
such that $(V\z)_\even:=V\oplus V^*$
and $(V\z)_\odd:=L^*$.
The assignment $\ell\mto j(\ell)\in T_RV\sset T_R(V\z)$
gives an $R$-super-bimodule imbedding $j\z : L\into T_R(V\z)$.
We  put $\ta:=T_R(V\z)/(L),$ where $(L)$ stands for the two-sided ideal
generated by the image of $j\z $.

\begin{lemma}\label{z} 
The complex $(\KK_\idot\ta)/L\lo$ is acyclic in positive
homological
degrees if and only if each of the two complexes $K_\idot{A}$ and
 $\KK_\idot{A}/L\lo$ is acyclic in positive
homological
degrees.
\end{lemma}
\begin{proof} An element
$m$ of  an  $R$ super-bimodule $M$ is called homogeneous
if $x\in M_\even,$ in which case we set $|x|=0$,
or  $x\in M_\odd,$ in which case we set $|x|=1.$
Given an integer $k\geq 1$,
let
$M^{\o k}\anti$ denote a quotient of 
$T^k_RM$ modulo all relations of the form
$$
m_1\o\ldots\o m_k=(-1)^q\cdot m_k\o m_1\ldots\o m_{k-1},\quad
r\cd m_1\o\ldots\o m_k=m_1\o\ldots\o m_k\cd r,
$$
where $r\in R$ and
$m_1,\ldots, m_k\in M$ are homogeneous elements, and
where we put $q:=$ $|m_k|\cdot(|m_1|+\ldots+|m_{k-1}|).$

Now, let
$B=R\oplus B_+$ and $C=R\oplus C_+$ be a pair of graded  augmented $R$ super-algebras
and let $B*_RC$ denote their free product over $R$. 
Clearly, we have 
\begin{eqnarray}\label{big}
&(B*_RC)\big/[B*_RC,B*_RC]\,&=\,
B/[B,B]\bigoplus C/[C,C]\bigoplus
(B_+\o_RC_+)\\
&&\bigoplus (B_+\o_RC_+)^{\o 2}\anti\bigoplus (B_+\o_RC_+)^{\o
3}\anti\bigoplus\ldots,
\nonumber
\end{eqnarray}
where in the super-case the commutators 
 are to be replaced by super-commutators.

We apply this in the special case  $B:=T_R(V^*\oplus L^*)$
and $C:=T_R(V\oplus L)$.
We have
$K_\idot\ta=T(V\oplus V^*\oplus L\oplus L^*)=B*_RC$,
where the factor $B$ sits in homological degree $0$. Therefore, from \eqref{big}
 we obtain 
$$\KK_\idot\ta=\frac{B}{[B,B]}\bigoplus \KK_\idot{A}\bigoplus  
\big(B_+\o_R(K_\idot{A})_+\big)\bigoplus
\big(B_+\o_R(K_\idot{A})_+\big)^{\o 2}\anti
\ldots,
$$
where $B$ is viewed as a DG super-algebra with zero differential
concentrated in homological degree zero. 

Hence, 
for  homology  we get
$$
H_\idot\big((\KK_\idot\ta)/L\lo\big)=\frac{B}{[B,B]}\bigoplus 
H_\idot\big(\KK_\idot{A}/L\lo\big)\bigoplus  
B_+\o_RH_\idot\big((K_\idot{A})_+\big)\bigoplus\ldots.
$$
We see that that the homology in the RHS is concentrated
in degree 0 if and only if 
each of the two complexes $\KK_\idot{A}/L\lo$ and
$K_\idot{A}$ has vanishing homology in degrees $>0.$
\end{proof}
\begin{remark}  The  above Lemma may be seen as an instance of 
cyclic Kunneth formula
for free products. 
\end{remark}

\begin{definition}\label{ascyc}
Let $M^\dd=\bigoplus_{p\in\Z,r\geq0}M^\dd_p[r],$
 be a collection of
bigraded complexes (one for each $\dd\in\Z_+^I$),
with differentials $d=d^\dd: M^\dd_\idot[r]\to M^\dd_{\hdot-1}[r]$.

We say that the complexes $(M^\dd, d)$ are
{\em asymptotically acyclic} in positive  homological  degrees
 as $\dd\to\infty$ if, for any
$r, N\geq 1$, there exists a dimension vector
$\dd=\dd(r,N)\succ N$
such that one has $H_k(\bigoplus_{s\leq r}M^\dd_\idot[s])=0$ for all
$k>0.$
\end{definition}

\begin{proposition}\label{dummy} For the graded super-algebra
$A=T_RV/(L)$ associated to a $(V,L)$-datum, the following are equivalent:

\vi The complexes $(\BK_\idot(\Rep_{\dd}\ta)^{G_{\dd}}, d_K),$
\textbf{of} $G_{\dd}$\textbf{-invariants},
are asymptotically acyclic in positive  homological  degrees as $\dd\to\infty$.

\vii The complexes $(\BK_\idot(\Rep_{\dd'} A), d_K)$
are asymptotically acyclic in positive  homological  degrees as $\dd'\to\infty$.
\end{proposition}

\begin{proof} 
The graded super-algebra $\ta$ is 
clearly isomorphic to
$A\,*_R\, T_R(V^*\oplus L^*),$ a free product (over $R$) of $A$,
and the tensor algebra $T_R(V^*\oplus L^*)$.

For any dimension vector $\dd$, we have
$\Rep_\dd T_R(V^*\oplus L^*)=\VV_\dd ^*\oplus\LL_\dd^*,$
where $\VV_\dd ^*\oplus\LL_\dd^*$ is the bigraded vector space dual to
$\VV_\dd\oplus\LL_\dd$.
 Therefore, 
 one has a canonical isomorphism
$$
\Rep_\dd\ta\cong\Rep_\dd A\times\Rep_\dd T_R(V^*\oplus L^*)=
\Rep_\dd A\times\VV_\dd ^*\times\LL_\dd^*\sset
\VV\times\VV_\dd ^*\times\LL_\dd^*,
$$
where all the sets above are viewed as {\em super}-schemes.
For  the corresponding  Koszul complexes,
this gives the following $G_\dd\times\C^\times$-equavariant
DG algebra isomorphism
\beq{ta}
\BK_\idot(\Rep_\dd\ta)\cong\BK_\idot(\Rep_\dd A\times\VV_\dd ^*\times\LL_\dd^*)
\cong\BK_\idot(\Rep_\dd A)
\otimes\big(\C[\VV_\dd ^*]\o\La\LL_\dd^*\big),
\eeq
where the rightmost tensor factor $\C[\VV_\dd ^*]\o\La\LL_\dd^*$
is a DG super-algebra equipped with
zero differential and concentrated in homological degree zero.

Now fix integers $r,N$ and assume that the complex
$(\bigoplus_{s\leq r}\BK_\idot(\Rep_\dd A)[s], d_K)$
is acyclic in positive  homological  degrees, for some
$\dd\succ N$.
We see from \eqref{ta}  that a similar statement holds for the
complex $\BK_\idot(\Rep_\dd\ta)$ as well. It follows that
the same holds for the subcomplex 
 $\BK_\idot(\Rep_\dd\ta)^{G_\dd}$, of $G_\dd$-invariants.
This proves
the implication (i)$\en\Rightarrow\en$(ii) of the Proposition.

To prove the opposite implication, we use  complete reducibility
of finite dimensional $G_\dd$-representations. Specifically,
 given a (finite dimensional)
rational  $G_\dd$-representation $U$ and an irreducible $G_\dd$-representation $M$,
let $U^M=\Hom_{G_\dd}(M,U)$ be the $M$-multiplicity space of $U$.
For any pair $U,W$, of finite dimensional
 $G_\dd$-representations,
by complete reducibility, one has a canonical  direct 
sum decomposition of the space of $G_\dd$-diagonal invariants:
$(U\o W)^{G_\dd}=\bigoplus_{M\in\text{Irr}(G_\dd)} U^M\o W^{M^*},$
where $M^*$ stands for the contragredient representation.

Now, the group $G_\dd$ acts diagonally on  the
tensor product in the RHS of \eqref{ta}. Hence, taking $G_\dd$-invariants
we obtain 
\begin{eqnarray}\label{kost}
&\BK_\idot(\Rep_\dd\ta)^{G_\dd}&=
\Big(\BK_\idot(\Rep_\dd A)\o\C[\VV_\dd
^*]\o\La\LL_\dd^*\Big)^{G_\dd}\\
&&=
\bigoplus_{M\in\text{Irr}(G_\dd)} 
\BK_\idot(\Rep_\dd A)^M\otimes\big(\C[\VV_\dd
^*]\o\La\LL_\dd^*\big)^{M^*}.
\nonumber
\end{eqnarray}

Observe also that, for each irreducible
$G_\dd$-representation $M$, one clearly has
\beq{kost2}
\BK_\idot(\Rep_\dd A)^M=
\big(\C[\VV_\dd]\o\La\LL_\dd\big)^{M}\neq 0
\en\Longleftrightarrow\en
\big(\C[\VV_\dd^*]\o\La\LL_\dd^*\big)^{M^*}\neq 0.
\eeq

To complete the proof, fix integers $r,N,$ and assume that for some
$\dd\succ N$ the complex
$(\bigoplus_{s\leq r}\BK_\idot(\Rep_\dd \ta)[s]^{G_\dd}, d_K)$
is acyclic in positive  homological  degrees.
We deduce from the direct sum decomposition
in \eqref{kost} and from \eqref{kost2} that,
increasing the values of $r$ and $\dd$ if necessary, 
one can ensure that a similar statement holds for each
$G_{\dd}$-isotypic component that occurs with nonzero multiplicity in the
complex $\bigoplus_{s\leq r}\BK_\idot(\Rep_\dd)[s]$.
The implication (ii)$\en\Rightarrow\en$(i)
follows.
\end{proof}

\subsec{Proof of Theorem \ref{converse}.} 
Given any
$(V,L)$-datum, we 
may view the tensor algebra $T_R(V\oplus L)$
as a free bigraded super-algebra such that the vector space
$V$ is even and is placed in homological degree 0,
resp. the vector space
$L$ is odd and is placed in homological degree 1
(cf. Proof of Proposition \ref{HC1}).

We apply the general construction of the
trace map $\TR$, cf. \eqref{trace},
to the free bigraded super-algebra
$K_\idot{A}=T_R(V\oplus L)$.
Since $
\Rep_\dd\big(T_R(V\oplus L)\big)\cong\VV_\dd \oplus \LL_\dd,$
 we obtain a bigraded super-algebra homomorphism
\beq{KTR}
\TR: \ \O(K_\idot{A})=\ssym\big((\KK_\idot{A})_+\big)
\too \big(\C[\VV_\dd ]\o\La^\hdot\LL_\dd^*\big)^{G_\dd}.
\eeq
It is straightforward to check that this map intertwines the Koszul
differentials on both sides, i.e., it is a morphism of DG algebras.

Next, recall the short exact sequence from \eqref{psi}
and the isomorphism in \eqref{BK_equiv}.
It is easy to see that  the following diagram commutes
\beq{lll}
\xymatrix{
L\lo\;\ar@{^{(}->}[rrrr]^<>(0.5){\eqref{KK}}\ar@{=}[d]^<>(0.5){\Id}
&&&&\;\K_1{A}\ar[d]^<>(0.5){\TR}\\
L\lo\;\ar@{^{(}->}[rr]^<>(0.5){\psi^\top}&&
\;(\LL_\dd^*)^{G_\dd}\;\ar@{^{(}->}[rr]
&&\;\big(\C[\VV_\dd]\o\La^1\LL_\dd^*\big)^{G_\dd}.
}
\eeq
Therefore,  the map $\TR$ descends to a well defined map
\beq{red}
\TS_\dd:\
\ssym\big((\KK_\idot{A})_+/L\lo\big)
\too \big(\C[\VV_\dd ]\o\La^\hdot(\LL_\dd^*/L\lo)\big)^{G_\dd}=
\BK_\idot(\Rep_\dd A)^{G_\dd}.
\eeq
(here and below, we will identify the space
$L\lo$ with its image $\psi^\top(L\lo)\sset\LL_\dd^*).$
 The map \eqref{red} is again a bigraded super-algebra homomorphism
compatible with the Koszul differentials on each side.

\begin{lemma}\label{TS}
The maps $\TS_\dd$ are
 asymptotically bijective, cf. Definition \ref{asbij}.
\end{lemma}
\proof Since any free algebra is RCI,  Proposition 
\ref{surjinj} implies that
the map $\TR$ in \eqref{KTR} is an
 asymptotically bijective algebra homomorphism.
By commutativity of diagram \eqref{lll},
this homomorphism sends the subspace $L\lo\sset
\O(K_\idot{A})$ isomorphically onto the
corresponding subspace $\psi^\top(L\lo)\sset \BK_\idot(\Rep_\dd A)^{G_\dd}.$
Therefore, the map $\TR$  yields an asymptotic
bijection between  the ideal
of the algebra $\O(K_\idot{A})$ generated by the subspace $L\lo$
 and   the ideal of the algebra
$\BK_\idot(\Rep_\dd A)^{G_\dd}$
 generated by the subspace $\psi^\top(L\lo)$.
The  induced  asymptotic bijection
$\O(K_\idot{A})/\O(K_\idot{A})\cdot L\lo\iso
\BK_\idot(\Rep_\dd A)^{G_\dd}/\BK_\idot(\Rep_\dd A)^{G_\dd}
\cdot\psi^\top(L\lo)$
is nothing but the map $\TS_\dd$, since
we have 
$$
(\C[\VV_\dd ]\o\La^\hdot\LL_\dd^*)^{G_\dd}\big/
(\C[\VV_\dd ]\o\La^\hdot\LL_\dd^*)^{G_\dd}\cdot L\lo
\cong\big(\C[\VV_\dd ]\o\La^\hdot(\LL_\dd^*/L\lo)\big)^{G_\dd}.\qquad\Box
$$
\vskip 4pt

\begin{proof}[Proof of Theorem \ref{converse}]
Assume first that the condition
of Theorem \ref{converse}(ii) holds. Thus, the complexes 
$(\BK_\idot(\Rep_{\dd'} A), d_K)$ are  asymptotically  acyclic 
in positive homological degrees as
$\dd'\to\infty$,
cf. Definition \ref{ascyc}.
Hence,
the complexes $(\BK_\idot(\Rep_{\dd}\ta)^{G_{\dd}}, d_K\big),$
of $G_{\dd}$-{\em invariants},
are asymptotically  acyclic in positive homological degrees  as $\dd\to\infty$,
by Proposition
\ref{dummy}. 
Applying Lemma \ref{TS}  to the algebra $\ta$, we deduce
that the complex $(\KK_\idot\ta)/L\lo$ is
acyclic in positive homological degrees. This implies,
by Lemma \ref{z}, that each of 
the complexes $K_\idot{A}$ and
$(\KK_\idot{A})/L\lo$ is acyclic in positive homological degrees.
Thus, $A$ is a NCCI, by Theorem \ref{Anick}(5). Now,
 the implication
(ii)$\en\Rightarrow\en$(i) of Theorem \ref{converse}
follows from  Lemma ~\ref{easy}.

Conversely, assume condition (i)  of Theorem \ref{converse} holds.
From   Lemma \ref{z}, we deduce that 
 the complex $(\KK_\idot\ta)/L\lo$ is
acyclic in positive homological degrees.
Hence,
the complexes 
$(\BK_\idot(\Rep_{\dd}\ta)^{G_{\dd}}, d_K),$
of $G_{\dd}$-{\em invariants},
are asymptotically  acyclic in positive homological degrees as
$\dd\to\infty$,
 by Lemma \ref{TS} applied to
the algebra $\ta$.
Therefore,  the complexes
$(\BK_\idot(\Rep_{\dd'} A),d_K)$ are asymptotically  acyclic  in positive
homological degrees
as $\dd'\to\infty$, by Proposition
\ref{dummy}. This proves
the implication
(i)$\en\Rightarrow\en$(ii).

Finally, the equations for Hilbert series given in Theorem
\ref{converse} 
follow from \eqref{ncci_h}.
\end{proof}

\section{Additional results and examples}\label{Ex}
\subsec{Twisted free products}\label{twfree} In this subsection,
we provide a few useful ways of
producing NCCI, resp. RCI algebras.

First, let $A_j=T_RV_j/(L_j)$, $j=1,2$, be a pair graded algebras associated with
$(V_i,L_i)$-data. 
The free product of $A_1$ and $A_2$ over $R$ is
 a graded algebra
$$A_1*_R A_2\cong T_R(V_1\oplus V_2)/(L_1\oplus L_2).$$

One easily finds
\beq{formula}
h(A_1*_R A_2)=\big(h(A_1)^{-1}+h(A_2)^{-1}-\bone\big)^{-1}.
\eeq

\begin{proposition} \label{sf}
\vi Let $A=T_RV/(L)$, and $B=T_RV/(L')$, where $L'\subset L$ is a
subbimodule. Then, if $A$ is RCI then $B$ is RCI. 

\vii  If $A_1*_R A_2$ is RCI then both $A_1$ and $A_2$
are RCI. 

\viii The converse to \vii holds if,
in Definition \ref {rci_def},
for each $N\gg0$  one can use the same
dimension $\mathbf d$  both for  $A_1$ and~$A_2$. 
\end{proposition}

\begin{proof} Observe that
if an algebraic map to a vector space is flat then so is the induced map 
to any quotient vector space. Part (i) follows.  Parts (ii)-(iii) 
are immediate from \eqref{dim}.
\end{proof}

Assume next that one has a $(V,L)$-data and
 an $(L,M)$-data, that is, there are 
 finite dimensional  $\mathbb   N$-graded $R$-bimodules
$V,L,M$ and  graded maps
\beq{jj} j_L: L\into T^+_RV,\quad\text{and}\quad
 j'_M:\xymatrix{
M\;\ar@{^{(}->}[rr]^<>(0.5){ j_M}&&T^+_RL\ar[rr]^<>(0.5){T( j_L)}&&
 T^+_RV},
\eeq
where $T( j_L)$ is the algebra morphism induced by the
$R$-bimodule morphism $ j_L.$
Thus, we may define a triple of graded $R$-algebras
\beq{BDA}
B:=T_RV/(\im j_L),\quad D:=T_RL/(\im j_M),\quad B\circ_R D:=T_RV/(\im j'_M).
\eeq

There is a slightly different but equivalent interpretation
of the algebra $B\circ_R D$. To explain it, consider two  $R$-bimodule maps
\begin{eqnarray*}
& j_{LM},\;j'_{LM}: L\oplus M \too T_R(V\oplus L),
\quad& j_{LM}(\ell\oplus m):=j_L(\ell)+ j_M(m),\\
&&j'_{LM}(\ell\oplus m):=  j_L(\ell)-\ell+ j_M(m).
\end{eqnarray*}
We observe that the algebra $B*_RD$, a  free product,
is a quotient of $T_R(V\oplus L)$
by the two-sided ideal
generated by the image of the map $ j_{LM}$.
 Observe further that 
the kernel
of the algebra homomorphism
$$
T_R(V\oplus L)\onto
T_RV/(\im j'_M)=B\circ_R D,
$$
induced by the first projection $V\oplus L \onto V,$
is a two-sided ideal in $T_R(V\oplus L)$ generated
by the image of   the map 
$j'_{LM}$.
Therefore, we get 
\beq{cst}
B*_RD= T_R(V\oplus L)/(\im j_{LM}),\quad\text{and}\quad
B\circ_R D= T_R(V\oplus L)/(\im j'_{LM}).
\eeq

\begin{proposition}\label{nccn} If $B$ is an NCCI algebra then,
for  the Hilbert series of $B\circ_R D$, we have
$$
h(B\circ_R D)=\big[\bone-h(D)\cd(h(V)-h(L))\big]\inv\cdot h(D).
$$ 
\end{proposition}

To prove the proposition, we introduce an increasing
filtration on the
algebra $B\circ_R D$ as follows. First define a grading
on $T_R(V\oplus L)$ by assigning  elements
of the vector space $V\sset V\oplus L$ their natural degrees
and placing the vector space $L\sset V\oplus L$
in degree zero. The increasing filtration on
 $T_R(V\oplus L)$ induced by this grading descends
to an increasing filtration on the quotient
algebra $B\circ_R D=
T_R(V\oplus L)/(\im j'_{LM})$. Let $\gr(B\circ_R D)$ denote the corresponding associated
graded algebra.

It is clear from the isomorphisms in \eqref{cst}
that the tautological imbedding
$V\oplus L\into T_R(V\oplus L)$
 induces a well defined and 
 surjective graded algebra  homomorphism
\beq{grbd}\xi:\ B*_RD =
T_R(V\oplus L)/(\im j_{LM})\onto {\rm gr}(B\circ_R D)=\gr\big(T_R(V\oplus L)/(\im j'_{LM})\big).
\eeq

\begin{lemma}\label{hilser1}  If $B$ is an NCCI algebra then the map
\eqref{grbd}
is an isomorphism.
\end{lemma}

\begin{proof}[Proof of Lemma]
Let $F=T_RV$, and to simplify notation write
$\otimes=\otimes_R, \,(-)*(-)=(-)*_R(-),$ etc.
The noncommutative
Koszul complex $K_\idot(B)$ for the algebra $B$ reads
$$
\ldots\to F\otimes L\otimes F\otimes L\otimes F\to F\otimes
L\otimes F\to F\onto B.
$$
Taking a free product of this complex with the algebra $D$
yields a complex of the form
\begin{equation}\label{S0}
\ldots\to (F*D)\otimes L\otimes (F*D)\otimes L\otimes (F*D)\to (F*D)\otimes
L\otimes (F*D)\to F*D\to B*D.
\end{equation}

The complex  $K_\idot(B)$ provides a resolution of $B$,
since $B$ is an NCCI algebra.
It follows that the complex in \eqref{S0} provides a resolution
for the algebra $B*D$.

On the other hand, using the presentation of the
algebra $B\circ D$ given in \eqref{cst}, one sees
that the terms of the corresponding Koszul complex
$K_\idot(B\circ D)$ can be written in the form
\begin{equation}\label{S}
\ldots\to (F*D)\otimes L\otimes (F*D)\otimes L\otimes (F*D)\to (F*D)\otimes
L\otimes (F*D)\to F*D\to B\circ D.
\end{equation}

Now, the above defined increasing filtration on $B\circ D$ 
gives rise to an increasing filtration on the DG algebra
$K_\idot(B\circ D)$. Furthermore, it is clear that
 the complex (\ref{S0}) is obtained by taking 
the associated graded of the complex (\ref{S}) with respect
to that filtration.  Hence, it follows that, for
an NCCI algebra $B$, the complex
(\ref{S}) is acyclic in positive homological degrees. 
This implies that $\xi$ is an isomorphism,
and the lemma is proved.
\end{proof}

\begin{proof}[Proof of Proposition
\ref{nccn}] Since $B$ is NCCI,  we
have $h(B)=(\bone-h(V)+h(L))^{-1}$. Further,
by Lemma \ref{hilser1}, the graded
$R$-algebras  $B\circ D$ and
$B*D$ have the same  Hilbert series.
The proposition now follows from  formula \eqref{formula}.
\end{proof}

\subsec{} We keep the notation of \eqref{jj}-\eqref{BDA}.
Also put $F:=T_RV,$ and  $A:=B\circ_R D$ and $A':=B*_R D.$
 We have the sub $R$-bimodule 
$L^\circ=L\cap [F,F]$. 

Let $Q:=L\cap [F,F]\cap (\im j_M)$
be
the kernel of the natural map $\eta: L^\circ\to D=T_RL/( j_M).$

\begin{proposition}\label{hilser2}
If $B$ is an asymptotic RCI algebra then one has 
$$
h({{\mathcal O}(B\circ_R D)};t)=\frac{h({{\mathcal O}(D)};t)}{
\lambda(Q)\cd\prod_{s\ge 1}\det\big(\bone +h(D;t^s)\cd[h(V;t^s)-h(L;t^s)]\big)}.
$$ 
\end{proposition}
 Let $\bar L^\circ$ be
the image of the  map $\eta: L^\circ\to D$.
We have a natural 
embedding $\gamma: \bar L^\circ\into [A,A]$. 
The image of $\gamma$ is contained in the degree zero part of $[A,A]$ 
under the filtartion on $A$, so it gives rise to an embedding 
$\gr\gamma: \bar L^\circ\into \gr([A,A])$.    

\begin{lemma}\label{asRCI} One has 
$\dis
\gr([A,A])=[A',A']\oplus (\gr\gamma)(\bar L^\circ).
$
\end{lemma}

\begin{proof} The direct sum of the terms of the complex
 \eqref{S}, resp. \eqref{S0}, gives a
 differential graded
algebra $K_\idot$, resp. $K'_\idot=K_\idot(B)*_RD$.
As we have explained earlier, the DG algebra  $K_\idot$ comes equipped with
the increasing filtration,
such that one has $\gr  K_\idot=K'_\idot.$
Therefore, for the corresponding super-commutator quotients,
$\KK_\idot:=K_\idot/[K_\idot,K_\idot]$ and
$(\KK_\idot)':=K'_\idot/[K'_\idot,K'_\idot]$, respectively, we deduce
 $\gr\KK_\idot=(\KK_\idot)'.$

 Since $B$ is an
asymptotic $RCI$ algebra, the homology of
the complex $\gr\KK_\idot=(\KK_\idot)'$ is
concentrated in degrees 0 and 1. Moreover, we have
$H_0((\KK_\idot)')=A'/[A',A']$ and
$H_1((\KK_\idot)')=L^\circ$.

The increasing filtration on $\KK_\idot$ gives rise
to a standard spectral sequence
with  $E_1$-term $H_\idot(\gr\KK_\idot)$ that converges to
the homology of $\KK_\idot$. It can be seen
from the analysis of the spectral
sequence  that this yields $H_1(K_\idot)=Q$. We also have
$H_0(K_\idot)=A/[A,A]$.
The lemma follows.      
\end{proof}

\begin{proof}[Proof of Proposition \ref{hilser2}]
By Proposition \ref{asRCI}, to compute the Hilbert series of ${\mathcal O}(A)$
it suffices to compute the 
Hilbert series of ${\mathcal O}(A')$. To compute the latter,
note that we have  
$$
\frac{A'_+}{[A',A']}=\frac{B_+}{[B,B]}\oplus \frac{D_+}{[D,D]}\oplus
\frac{T_R(B_+\otimes_R D_+)_+}{[T_R(B_+\otimes_R D_+),T_R(B_+\otimes_R D_+)]}.
$$
This implies that 
$$
h({{\mathcal O}(A')};t)=
\frac{h({{\mathcal O}(B)};t)\cdot
h({{\mathcal O}(D)};t)}{
\prod_{s\ge 1}\det\big(h(B;t^s)+h(D;t^s)-h(B;t^s)\cdot h(D;t^s)\big)}.
$$ 

Using that $B$ is an asymptotic RCI algebra, the last equation can be rewritten
as 
$$
h({{\mathcal O}(A')};t)=\frac{1}{\lambda(Q)}\frac{
h({{\mathcal O}(D)};t)}{ 
\prod_{s\ge 1}\det \big(\bone+h(D;t^s)\cdot[h(V;t^s)-h(L;t^s)]\big)},
$$ 
and the proposition follows.
\end{proof}

\subsec{Examples of RCI algebras.} We begin with the example
of the algebra of $q$-polynomials. Thus, we fix
 $q\in \mathbb \C^\times$, and
let $\mathbb C_q[x,y]$ be  the  algebra 
 with two generators
$x,y$ of degree $1$ 
and defining relation $xy-qyx=0$.

\begin{proposition}\label{k2} If  $q$ is not a root of unity,
then
the  algebra  $\mathbb C_q[x,y]$ is RCI.
\end{proposition}

\begin{remark}\vi If $q$ is a root of unity then the $q$-polynomial algebra
$\mathbb C_q[x,y]$ is not a RCI algebra. Namely, if $n$ is the order
of $q$ then $\mathbb C_q[x,y]$ has a 2-parameter family of
equivalence classes of irreducible $n$-dimensional
representations, which gives an $n^2+1$-dimensional family of
representations. Taking direct sums of such representations and 
1-dimensional representations, we see that $\Rep_d(A)$ has
dimension  bigger than expected, for any $d\ge n$.

\vii Formulas \eqref{AA}(i)-(ii) for 
$\mathbb C_q[x,y]$ are easy to check directly. 
If $q$ is a root of unity, then it is easy to see that 
(i) still holds but (ii) fails (which is another way to see that
in this case $\mathbb C_q[x,y]$ is not RCI).
\end{remark}

\begin{proof}\footnote{This very short proof was explained to us
by G. Lusztig}
We need to show that the variety $R_n$ of $n$-dimensional representations
of the algebra $\mathbb C_q[x,y]$ has dimension $n^2$. 
To do so, consider the subvariety $R_n^0$ of $R_n$ 
consisting of representations on which $x$ is a nilpotent matrix. 
Since $x$ in this case is conjugate to $qx$, 
the dimension of the space of solutions $y$ of the equation 
$xy=qyx$ is the same as the dimension of the centralizer of $x$. 
Therefore, for each nilpotent conjugacy class $C$ (corresponding
to the matrix $x$), 
the dimension of the set $\{(x,y)\mid xy=qyx,\;x\in C\}$
 equals $n^2$. Thus, we get $\dim R_n^0=n^2$. 

Any representation $V\in R_n,$ of $\mathbb C_q[x,y]$,
has a canonical direct sum decomposition into
Fitting components $V=V_1\oplus V_2$, such that the endomorphism
$x|_{V_1}$ is nilpotent and  the endomorphism
$x|_{V_2}$ is invertible (and then $y$ is nilpotent because $q$ is not a
root of unity). Let $\dim V_1=n_1,$ and $\dim V_2=n_2$. 
Then there are $(n_i)^2$ parameters for $V_i$, and an additional 
amount of $2n_1n_2$ parameters for the choice of the vector space
decomposition $V=V_1\oplus V_2$. Thus the total dimension is
$(n_1+n_2)^2=n^2$, and we are done.  
\end{proof}

Here is a more interesting example of a RCI algebra
where formulas \eqref{AA} are not
trivial.

\begin{example} Let $m\geq 1$, and fix an $m$-tuple
$q_1,\ldots,q_m\in\C^\times.$
Let $A$ be the algebra  with $m$ generators 
$x_i$, $i=1,\ldots,m$, and one relation 
$$
{\rm ad}_{q_1}(x_1){\rm ad}_{q_2}(x_2)\ldots{\rm ad}_{q_{m-1}}(x_{m-1})x_m=0,
$$
where ${\rm ad}_q(x)y:=xy-qyx$, and $q_i\in \mathbb \C^\times$ are not
roots of 1. We claim that it is a RCI algebra. 
Indeed, we can rewrite the relations as 
\begin{eqnarray*}
&x_{m-1}x_m-q_{m-1}x_mx_{m-1}=z_{m-1},\quad
&x_{m-2}z_{m-1}-q_{m-2}z_{m-1}x_{m-2}=z_{m-3},\;\ldots,\\
&&x_1z_2-q_1z_2x_1=0.
\end{eqnarray*}

We need to show that the representation space of this algebra has
dimension $n^2(m-1)$. We assign degree 1 to all
variables, $x_i,z_j$. The  corresponding top degree
homogeneous terms  of the above equations read
\begin{eqnarray*}
&x_{m-1}x_m-q_{m-1}x_mx_{m-1}=0,\quad
&x_{m-2}z_{m-1}-q_{m-2}z_{m-1}x_{m-2}=0,\;\ldots,\\
&&x_1z_2-q_1z_2x_1=0.
\end{eqnarray*}
we find that this defines a complete intersection by 
Proposition \ref{k2} (so the dimension of the space of solutions is the
expected number $n^2(m-1)$). Therefore, the original equations 
also have the space of solutions of this dimension, as desired. 
\end{example}

\subsec{Examples of asymptotic RCI algebras.}\label{asex}
We begin with some combinatorial preliminaries. 
Let $e$ and $f$ be two (possibly equal) words in the alphabet $v_1,\ldots,v_n$. 
We will say that $e,f$ are non-overlapping if 
none of them is a proper subword of the other, and 
none of them begins the way the other one ends. 
It is clear that $e,f$ are non-overlapping if and only if
any word $w$ that contains both $e$ and $f$ (in different ways
if $e=f$) is of the form $w_1ew_2fw_3$ or $w_1fw_2ew_3$.

Let  $V$ be a vector space with basis $v_1,\ldots,v_n$.
We have a free $\C$-algebra 
$TV=
\C\langle v_1,\ldots,v_n\rangle,$
where $v_k$ is assigned some grade degree $m_k$.
Given homogeneous elements
$e_1,\ldots,$
$e_p\in TV$,  we put
$A=TV/(L)=
\C\langle v_1,\ldots,v_n\rangle/(e_1,\ldots,e_p).$

\begin{theorem}\label{kcexact}
Suppose that the defining relations $e_1,\ldots,e_p$ of $A$ 
are monomial (i.e. each is given by a single word)
and pairwise non-overlapping. Then $A$ is an  asymptotic ~RCI.
\end{theorem}

 Note that the condition that $e_1,\ldots,e_p$ be pairwise
non-overlapping includes the stipulation that 
each $e_i$ is non-overlapping with itself. 
Anick uses the name `{\em strongly free}' for such a 
 pairwise
non-overlapping set  $e_1,\ldots,e_p$, see \cite{An}.
\begin{proof} First, one shows by a standard
inclusion-exclusion argument that Anick's complex
$C_\idot{A}$ is a resolution of $A$. Thus, $A$ is NCCI.

Next, observe that in our case, we have
$L\lo=0$ (this follows from the fact that the relations are
non-overlapping).
Hence, by Theorem \ref{converse},
we must only prove that the complex $\KK_\idot{A}$ has vanishing
homology in all positive degrees.

\begin{remark} The vanishing of the homology groups of
$\KK_\idot{A}$ of degrees $>1$ is a consequence of Proposition
\ref{HC1}, but it may also be proved directly as follows.
\end{remark}

For brevity, write
$Q_i=\KK_i{A}$, $i\ge 0$.
We have a natural decomposition
$Q_\idot=\oplus_w Q_\idot(w)$, where the summation
is taken over cyclic words $w$, and $Q_i(w)$ for $i>0$
is the span of all elements $w_1\otimes e_{j_1}\otimes\ldots\otimes
w_i\otimes e_{j_i}$ such
that $w_1e_{j_1}\ldots w_ie_{j_i}=w$, while $Q_0(w)$ is a 1-dimensional
vector space
spanned by $w$.

Let $w$ be a cyclic word. Since the relations are
non-overlapping, the word $w$ has a unique (up to cyclic
permutations) maximal length representation as $w=z_1\ldots z_m$,
where $z_i:=w_ie_{j_i}$ (in this case, $w$ contains exactly $m$
distinct subwords of the form $e_i$). Let $D_\idot$ be the 2-step complex
with $D_0=D_1=\mathbb C$, and $d: D_1\to D_0$ the identity
map. It is easy to check  that if $w$ contains at least one
monomial $e_i$, and $\Gamma=\mathbb Z/n\mathbb Z$ is the cyclic
symmetry group of $w$ (clearly, $m$ is divisible by
$n$), then the complex $Q_\idot(w)$ is isomorphic to
$((D_\idot)^{\otimes m})^\Gamma$ and therefore is exact. On the other
hand, if $w$ does not contain $e_i$, then $Q_\idot(w)$ is
concentrated in degree 0. This implies that the complex $Q_\idot$
is exact in positive degrees, as desired. Theorem
\ref{kcexact} is proved.
\end{proof}

\begin{definition} Fix an alphabet $v_1,\ldots,v_n$
and let  $(m_1,\ldots,m_n; r_1,\ldots,r_p)$ be a collection of nonnegative
degrees. For each $i=1,\ldots,n,$ put $\deg v_i=m_i$.
We say that the above collection of degrees
 is {\em admissible} if 
there exist pairwise non-overlapping words $w_1,\ldots w_p$,
in our alphabet, such that $\deg w_j=r_j, \, j=1,\ldots p.$
\end{definition}

Let  $v_1,\ldots,v_n$ be indeterminates  of degrees $m_1,\ldots,m_n$. 
One says that a certain statement holds for a {\em  Weil generic}
$p$-tuple of homogeneous elements
$e_1,..,e_p\in \mathbb C\langle v_1,\ldots,v_n\rangle,$ of degrees 
$r_1,\ldots,r_p$ if it holds for 
$p$-tuples which belong to an at most countable intersection of
nonempty Zariski open sets in the direct sum of the 
homogeneous components of the algebra
$\mathbb C\langle v_1,\ldots,v_n\rangle$
of degrees $r_1,\ldots,r_p$. 

\begin{theorem}\label{generic} Let $(m_1,\ldots,m_n; r_1,\ldots,r_p)$ 
be an admissible collection of
degrees and put $\deg v_j=m_j,\,j=1,\ldots,n.$ Then, for  a
 Weil generic set of  homogeneous relations
 $e_1,..,e_p\in \mathbb C\langle
v_1,\ldots,v_n\rangle,$  of degrees 
$r_1,\ldots,r_p$ respectively, 
the algebra 
$\mathbb C\langle v_1,\ldots,v_n\rangle/(e_1,\ldots,e_p)$ is an
asymptotic RCI.
\end{theorem}

\begin{proof} For each $j\geq 0,$ the property 
of the complex $\KK_\idot{A}$ to have vanishing higher homology groups
in degree $j$ is clearly an open condition. Thus, Theorem \ref{generic}
is an  immediate consequence of 
 Theorem \ref{kcexact}.
\end{proof}

We refer to \cite{An3} and \cite[Theorem 2.25]{Pi} for related results on 
NCCI algebras with generic relations.

\begin{example} The collection $(1,1;r)$ is admissible for any
$r>1$. Indeed, if $x,y$ are generators then the word 
$x^{r-1}y$ is not self-overlapping. Thus, an algebra $A$ with 
two generators of degree 1 and one Weil generic relation of degree $r$ 
is an asymptotic RCI, by Theorem \ref{generic}.
However, such an algebra is {\em not} RCI for any $r\ge 4$.
Indeed, in this case there are 4-dimensional representations 
where the generators $x,y$ act by any 
strictly upper triangular matrices in some
common basis. The dimension of the set of such representations
is 18. Thus the dimension of the space of $N\ge 4$ dimensional 
representations which are isomorphic to the sum of 
one of the above 4-dimensional 
representations and a bunch of 1-dimensional representations 
is $N^2+2$, which is bigger than $N^2$, required by the RCI
property.

Now consider collections $(1,1;r,s)$, $r,s>1$. If $r=2,3$,
then it is easy to see that such a collection is never
admissible. If $r=4$, then the smallest $s$ for which 
this collection is admissible is $s=5$; the corresponding words
are $x^2y^2$ and $xyxy^2$.
For randomly chosen relations of degrees
$(1,1;4,4)$, the  corresponding algebra is not NCCI.

\end{example}
\section{Preprojective algebras and Quiver varieties}\label{yang}

\subsec{Hilbert  series  for preprojective
algebras}\label{hilbpi}
An important example, which was the main motivation for
this study, is that of preprojective algebras.

Recall the setting of \S\ref{uninak}.
Thus, $\Pi$ is the preprojective algebra of a connected quiver $Q$
with vertex set $I$. We have $\Pi=T_RV/(L)$, where $V$
is a $\C$-vector space with basis formed
by the edges of the quiver $\overline{Q}$,
and $L$ is an $R$-subbimodule in $V\o_RV$
generated by the element
$\sum_{a\in Q}[a,a^*]$.
Therefore, we get $h(V; t)=t\cdot\bc$, and $h(L; t)=t^2\cdot\bone.$ 
In this case, one  can rewrite
the $\ze$-function in \eqref{zeta} in terms of
 the  adjacency matrix $\bc$ of the double $\QQ$, as follows
\beq{zeq}
\ze(Q;t):=\ze(V,L;t)=\prod_{s\ge 1}\frac{1}{\det (\bone-t^s\cd\bc+t^{2s}\cd\bone)}.
\eeq
The vector space $L\lo$ is spanned by a single
degree 2 element $\sum_{a\in Q}[a,a^*]$, hence,  we have
$\lambda({L\lo})=1-t^2$.

According
to \cite{CB} (see explanations in \cite{CBEG},
section 11), we know that   $\Pi$ is an RCI algebra,
provided the quiver $Q$ is  neither
Dynkin nor extended Dynkin.
Thus, Theorem \ref{nak} follows from Theorem
\ref{maint}.

\begin{example} Fix an integer $g\geq 1$. For any $n\geq 1$,
let  $A_{g,n}$ be the algebra with generators $x_i,y_i,\,
i=1,\ldots,g$, and one defining relation 
$(\sum [x_i,y_i])^n=0$.

In the  special case $n=1$, the corresponding
 algebra  $A_{g,1}$ is nothing but
the  preprojective algebra
of  a quiver with one vertex and $g$
edge-loops. This preprojective algebra
may be thought of as an `additive' (i.e. Lie algebra) analogue of the group algebra of
the fundamental group of a genus $g$ Riemann surface.

\begin{proposition}\label{riemsur}
Let $g>1$. Then, for $n=1$, we have
$$
h(A_{g,1})=\frac{1}{1-2g\cd t+t^2}, \quad\text{and}\quad
h(\O(\Pi))=\frac{1}{1-t^2}\prod_{s=1}^\infty \frac{1}{1-2g\cd t^s+t^{2s}}.
$$

For any $n\geq 2,$  we have
\begin{eqnarray*}
&h(A_{g,n})&=\frac{1-t^{2n}}{1-2gt+2gt^{2n+1}-t^{2n+2}}, \\
&h({{\mathcal O}(A_{g,n})})&=\prod_{s=1}^\infty 
\frac{1-t^{2(n-1+s)}}{1-2gt^s+2gt^{(2n+1)s}-t^{(2n+2)s})}.
\end{eqnarray*}
\end{proposition}

\begin{proof} In the case $n=1$, the formulas follow
  from Theorem \ref{nak}.

 If $n\geq 2,$ we apply the results of \S\ref{twfree}.
Specifically, in the
notation of that section,
we take $B=A_{g,1}$ and $D=\mathbb   C[z]/(z^n)$, and also $Q=0$. 
Thus, one has
$$h(D)=\frac{1-t^{2n}}{1-t^2}, \quad\text{and}\quad
h({{\mathcal O}(D)})=\prod_{i=1}^{n-1}\frac{1}{1-t^{2i}}. 
$$

The algebra $A_{g,1}$ is an asymptotic RCI and we have
$A_{g,n}=A_{g,1}\circ_R D$.
Thus, from the formulas for $n=1$,
applying  Proposition \ref{nccn} and Proposition \ref{hilser2},
we get the required formulas for $n>1$.
\end{proof}

We note that the formula of Proposition \ref{riemsur} for
 $h(\O(A_{g,1}))$ fails 
if $g=1$. 
Theorem \ref{nak} cannot be applied to that case,
since the quiver with one vertex and one edge-loop is an affine Dynkin
quiver. In this case, the representation scheme $\Rep_d(A_{g,1})$ 
is the commuting scheme, consisting of pairs of matrices $X,Y$
of size $d$ such that $XY-YX=0$. This scheme has dimension
$d^2+d$ and is not a complete intersection  unless
$d=1$. 

Note also that, for $n>1$, the algebra  $A_{g,n}$ (as
well as the algebra $D=\mathbb   C[z]/(z^n)$) is not NCCI and it has
infinite cohomological dimension. \hfill$\lozenge$
\end{example}

Next, let $Q$ be an arbitrary quiver  and
 let $J\subset I$ be a subset of the set of 
vertices of $Q$. Following \cite{EE}, define the {\it partial preprojective algebra}
$A=\Pi_{Q,J}$ in the same way as $\Pi_Q$ except that the defining
relation has the form 
$$
(1-p_J)\sum [a,a^*](1-p_J)=0, 
$$
where $p_J$ is the sum of the idempotents corresponding to the
vertices from $J$. 

\begin{proposition}\label{par}
If $J\ne \emptyset$ then, for any connected quiver $Q$,
 the partial preprojective algebra $\Pi_{Q,J}$ 
is RCI.
\end{proposition}

\begin{proof} If $Q$ has only one vertex, the statement is clear,
so we will assume that $Q$ has more than one vertex. 
Let $j\in J$, and let $Q'$ be the quiver $Q$ 
with an additional self-loop $\gamma$ at the vertex $j$. 
Then $\Pi_{Q',J}=\Pi_{Q,J}*_R\Pi_{\gamma,J}$, where $\gamma$
stands for the quiver with vertex set $I$ and unique edge
$\gamma$. The quiver $Q'$ is not Dynkin or affine Dynkin. 
Therefore, by Theorem \ref{nak}, $\Pi_{Q'}$ is a RCI
algebra, and hence by Lemma \ref{sf}, (i), $\Pi_{Q',J}$ is a RCI algebra. 
Hence by Lemma \ref{sf}, (ii), so is $\Pi_{Q,J}$, as desired. 
\end{proof}

Thus we obtain the following result.

\begin{proposition}\label{prepro1} Let
$Q$ be a connected quiver with vertex set $I$
and let $J\sset I$ be a {\em nonempty} subset.
Let $\bc$ be the adjacency matrix of
$\overline{Q}$ and $\bone_{I\sminus J}$ be the diagonal 
matrix whose principal diagonal is the characteristic
function of the set $I\sminus J$.  Then, we have 
$$
h(\Pi_{Q,J})=(\bone-\bc\cd t+t^2\cd \bone_{I\sminus J})\inv, 
\quad\text{and}\quad
h({\O(\Pi_{Q,J})})=\prod_{s=1}^\infty
\frac{1}{\det(\bone-t^s\cd \bc+t^{2s}\cd \bone_{I\sminus J})}.
$$
\end{proposition}

\begin{remark} \vi The  formula for $h(\Pi_{Q,J})$
was proved in \cite{EE} by another method.

\vii Note that there is no factor $(1-t^2)^{-1}$ in the
  formula for $h(\O(\Pi_{Q,J}))$
since $L\lo=0$ in this case. 
\end{remark}

\subsec{Asymptotic Hilbert series for Nakajima quiver varieties.}\label{mdw}
Fix a  pair of dimension vectors
$\dd=\{d_i\}_{i\in I}$ and $\ww=\{w_i\}_{i\in I}$.
For each $i\in I$, we
put 
$$\PP_i=\Hom(\C^{d_i},\C^{w_i}),\quad
\PP_{\dd,\ww}=\bigoplus_{i\in I} \PP_i,\quad
\BQ_i=\Hom(\C^{w_i},\C^{d_i}),\quad
\BQ_{\ww,\dd}=\bigoplus_{i\in I} \BQ_i.
$$

We write elements of $\PP,$ resp. of $\BQ,$ as
$ \bp=\{p_i\in \PP_i\}$, resp. $\bq=\{q_i\in\BQ_i\}.$ 
Thus, the composite $q_i\ccirc p_i$ is 
a linear map $\C^{d_i}\to\C^{d_i}$,
that is, an element of $\EE_{ii}=\Hom(\C^{d_i},\C^{d_i}).$

Now, let $Q$ be a connected quiver.
For any edge $a\in\QQ$, we put $\epsilon_a:=1$ if $a\in Q$,
resp.  $\epsilon_a:=-1$ if $a\in \QQ\sminus Q.$
Let $V$ be an $R$-bimodule
with basis formed by the edges of $\QQ$.
We have $\Rep_\dd(\QQ)=\Hom_{\bimod{R}}(V,\EE_\dd)=\VV_\dd.$

Nakajima considers a moment map
\begin{eqnarray*}
&\mu_{\dd,\ww}:\
\VV_\dd
\times\PP_{\dd,\ww}\times\BQ_{\ww,\dd}\to \oplus_{i\in I} \EE_{ii},\quad
(\rho,\bp,\bq)\mto \oplus_{i\in I}\,\mu_i(\rho,\bp,\bq),\\
&\mu_i(\rho,\bp,\bq):=
\left(
\sum\nolimits_{\{a\in\QQ\mid \text{tail}(a)=i\}} \epsilon_a\cd \rho(a^*)\cd
\rho(a)\right)\,-\,q_i\ccirc
p_i\,\,\in \EE_{ii},\quad\forall i\in I,\,\rho\in\Rep_\dd(\QQ).
\end{eqnarray*}

The affine quiver variety with parameters $\dd$ and $\ww$ is defined
as a categorical quotient
$\FM(\dd,\ww):=\mu_{\dd,\ww}\inv(0)//G_\dd.$
There is an action of the group $G_\ww=\prod GL(\C^{w_i})$
on $\FM(\dd,\ww)$ induced by the natural $G_\ww$-action on $\C^\ww.$
There is also
a $\C^\times$-action
on  $\FM(\dd,\ww)$ arising from 
the standard  $\C^\times$-action on
the vector space $\VV_\dd
\times\PP_{\dd,\ww}\times\BQ_{\ww,\dd}$, by dilation.
Thus, we get a weight  grading 
$\C[\FM(\dd,\ww)]=\bigoplus_{r\geq 0}\C[\FM(\dd,\ww)][r],$
such that, for each $r\geq 0$,
the weight component $\C[\FM(\dd,\ww)][r]$
 is a finite dimensional vector
space equipped with  $G_\ww$-action.

For any element $\bg=(g_1,\ldots,g_r)\in G_\ww=\prod GL(\C^{w_i})$
we form a generating function
$$\Phi\big(\C[\FM(\dd,\ww)];\bg,t\big)=
\sum\nolimits_{r\geq 0} t^r\cd\Tr\big(\bg|_{\C[\FM(\dd,\ww)][r]}\big)\in\C[[t]].
$$

Further, for each pair $i,j\in I$, we define
$\WW_{ij}=\Hom(\C^{w_i},\C^{w_j})=(\C^{w_i})^*\o \C^{w_j}
$.
Given $g_i\in GL(w_i)$ and $ g_j\in GL(w_j),$
we have a map $g_{ij}=g_i^\vee\o g_j: (\C^{w_i})^*\o \C^{w_j}\to
(\C^{w_i})^*\o \C^{w_j}.$ Thus, $g_{ij}\in \End\WW_{ij},$ and
we get a sequence of  polynomials
$$\det(\bone_{\WW_{ij}}-t^m\cd g_{ij})\in\C[t],\quad
m=0,1,\ldots.
$$

Let
$\phi_k\in\C[x],\,k=0,1,\ldots,$ be
 a  sequence 
of polynomials (essentially, Chebyshev polynomials
of the second kind) in an indeterminate $x$
 defined as the coefficients in
 the  expansion of the generating function
\beq{phi}\sum_{k=0}^\infty
\phi_k(x)t^k=\frac{1}{1-tx+t^2}.
\eeq
We set $\phi_k=0$ for $k$ negative.
Then, for each $k\in\ZZ$, we have an $I\times I$-matrix
$\|\phi(\bc)_{ij}\|$.
Further, recall the notation
$\ze(Q,t)$ from \eqref{zeq}.

An argument similar to that used in the proof of \cite{CBEG}, Lemma 11.1.2,
implies that, for any  
non-Dynkin quiver and any given  dimension vector $\mathbf w\ne 0$,
there exists  a sequence of dimension vectors $\mathbf d\to\infty$ 
such that the corresponding varieties $\mu_{\dd,\ww}\inv(0)$
are all complete intersections. Thus, repeating the  arguments 
of \S\ref{hilb_rci}
one obtains, via a matrix integral calculation similar to one
carried out in \S\ref{asymptotic},  the following version of 
Proposition \ref{hrep} for quiver varieties.

\begin{proposition}\label{nakprop} Fix a connected non-Dynkin quiver
$Q$ and a dimension vector $\ww\ne 0.$ For any $\bg\in G_\ww$ and $r\geq 0,$
the sequence of $r$-th coefficients in the formal
power series $\Phi\big(\C[\FM(\dd,\ww)]\big)$
 stabilizes as $\dd\to\infty,$ and for the
corresponding limit we have
$$
\underset{^{\dd\to\infty}}\lim \Phi\big(\C[\FM(\dd,\ww)];\bg,t\big)=
\ze(Q,t)\cdot
\prod_{i,j}\prod_{k\ge 0}\frac{1}{\det(\bone_{\WW_{ij}}-
t^{k+2}\cd g_{ij})^{\phi_k(\bc)_{ij}}}.\qquad\Box
$$ 
\end{proposition}

Next, we set 
$\langle\ww,\phi_k(\bc)\ww\rangle:=\sum_{i,j\in I}
\phi_k(\bc)_{ij}\cdot w_i\cdot w_j.$
Taking $\bg=1$ in the formula of the above Proposition,
we deduce
\begin{corollary} With the notation and assumptions  of Proposition
\ref{nakprop}, for
Hilbert series of quiver varieties
$\FM(\dd,\ww)$ we have 
$$
\underset{^{\dd\to\infty}}\lim h\big(\C[\FM(\dd,\ww)];t\big)=
\prod_{s\ge 1}\,\frac{1}{\det (\bone-t^s\cd\bc+t^{2s}\cd\bone)}\cdot
\prod_{k\ge 0}\,\frac{1}{(1-t^{k+2})^{\langle\ww,\phi_k(\bc)\ww\rangle}}.
\qquad\Box
$$
\end{corollary}

\begin{remark}
It was kindly pointed out to us by H. Nakajima that,
in \cite{NY}, the authors study  a function similar
 to $\Phi\big(\C[\FM(\dd,\ww)];\bg,t\big)$.
The problem solved in \cite{NY} is, however, quite different.
Specifically, the setup of  \cite{NY}  corresponds to the special
case where $Q$ is  the Jordan quiver. 
In that special case, Nakajima and Yoshioka analyze
asymptotic  behavior with respect to the variable
$t$ (corresponding, in the notation of
 \cite{NY}, to $e^{\varepsilon_1+\varepsilon_2}$)
for {\em fixed} dimension vector $\dd$.
The resulting  function (for  finite $\dd$)
is much more complicated than the one given in Proposition
\ref{nakprop}; it is
 a $K$-theoretic version
of the Seiberg-Witten prepotential, see 
\cite{NY} for details.
\end{remark}

\begin{remark} 
Fix a total ordering on $I$, the set of vertices of $Q$, and put
$GS_Q:=\prod_{i<j}GL(\bc_{ij})\times \prod_{i\in I}
Sp(\bc_{ii})$.
The
 group $GS_Q$ acts naturally on the vector space
spanned by the edges of $\QQ$ by linear transformations.
This gives a $GS_Q$-action on the path algebra of $\QQ$
and on the algebra $\Pi(Q)$,
by algebra automorphisms.
Therefore, for any pair $(\dd,\ww)$, of dimension vectors,
one has a linear  $GS_Q$-action on the symplectic vector space
$\VV_\dd\times\PP_{\dd,\ww}\times\BQ_{\ww,\dd}$
by symplectic automorphisms. This action commutes
both with the $G_\dd$-action and $G_\ww$-action,
as well as with the $\C^\times$-action.
Thus, the quiver variety $\FM(\dd,\ww)$ 
acquires the structure of a $G_\ww\times\C^\times\times GS_Q$-variety.

Now, let $\RR(GS_Q)$ denote the representation ring
of finite dimensional rational $GS_Q$-modules.
One can introduce
$GS_Q$-equivariant refinements of the formal
 series $h(\Pi;t),$
 $\Phi\big(\C[\FM(\dd,\ww)];\bg,t\big),$ etc.,
as appropriately defined elements of $\RR(GS_Q)[[t]]$,
the ring of formal power series with  $\RR(GS_Q)$-coefficients.
The formulas of  Theorem \ref{nak} and of Proposition
\ref{nakprop} still apply in this  refined $GS_Q$-equivariant
framework,
provided the symbol
$\bc_{ij}$ is understood as an element
of $\RR(GS_Q)$ corresponding to the class
of the fundamental  vector
representation of the group $GL(\bc_{ij})$ if $i<j$, 
resp.,  to  the class
of the fundamental covector
representation
of the group $GL(\bc_{ji})$ if $i>j$, and also as the class
of the vector representation of the group
$Sp(\bc_{ii})$
if $i=j$ (this way, $\bc$ becomes an
 $\RR(GS_Q)$-valued matrix).
\end{remark}

\subsec{The case of extended Dynkin quivers.}\label{extended}
Let $Q$ be an extended Dynkin quiver, and write
$o\in I$ for an extending vertex. 
We let $e_o$ denote the idempotent in the
preprojective algebra of $Q$ corresponding
to the trivial path at $o$.

Let $\G\sset SL_2(\C)$ be the finite subgroup
corresponding to $Q$ via the McKay correspondence.
Thus, the extending vertex corresponds to the trivial representation of 
$\G$.
We write $M$ for the 2-dimensional tautological representation of 
$\G$ and let  $\G\ltimes\sym M$ be the corresponding smash-product
algebra.

According to \cite{CBH}, one has
\begin{proposition}\label{cbh}\vi The  preprojective
algebra $\Pi$  is  Morita equivalent to  $\G\ltimes\sym M$.

\vii There is a graded algebra isomorphism
$\dis (\sym M)^\G\cong e_o\Pi e_o.$\qed
\end{proposition}

Part (i) of the Proposition implies, in particular,
that  Anick's complex is a resolution for the
preprojective algebra $\Pi$, hence, 
 $\Pi$ is  NCCI.
Therefore, Theorem \ref{Anick}(2) and the second formula in \eqref{ncci_h} yield
\beq{anick_pi}
h(\Pi)=(\bone-t\cd\bc+t^2\cd\bone)\inv,\quad\text{resp.},
\quad
 h\big(\O(\Pi)\big)=\ze(Q,t)\cdot h\big(\sym {HH}_2(\Pi)\big).
\eeq

Recall Chebyshev polynomials
of the {\em first kind}, $\varphi_k$, introduced in \eqref{varphi},
and  Chebyshev polynomials
of the {\em second kind}, $\phi_k$, introduced in \eqref{phi}.
These polynomials are related by the formula
$\varphi_k=\phi_k-\phi_{k-2},\,k=0,1,\ldots.$

We have the following well known result, cf. \cite{Su}. 

\begin{corollary}\label{b}
The Hilbert series of the graded algebra $(\sym M)^\G$ is given by
the formula
$$h\big((\sym M)^\G\big)=\sum\nolimits_{k\geq 0} \phi_k(\bc)_{oo}\cd t^k.
$$
\end{corollary}
\begin{proof} We use
the first formula in \eqref{anick_pi}, which  is a matrix equation.
Hence, for the
matrix element on each side of that equation corresponding to the
extending vertex,  we get
$$
h(e_o\Pi e_o)=h(\Pi)_{oo}=\big((\bone-t\cd\bc+t^2\bone)\inv\big)_{oo}=
\sum\nolimits_{k\geq 0} \phi_k(\bc)_{oo}\cd t^k,
$$
where the last equality is the definition of the polynomials $\phi_k$
given in \eqref{phi}.
The result now follows from  the isomorphism of Proposition
\ref{cbh}(ii).
\end{proof}

\begin{proof}[Proof of identity \eqref{curious}.]
We will be interested in Hochschild cohomology
of the  smash-product $\G\ltimes\sym M$. These 
are well known and easy to find.
To describe the answer,  put $\SS:=\G\sminus\{1\}$, and let
$\C[\SS]^\G$ denote the vector space of
class functions on the finite set $\SS$.
Then, one has, see \cite{CBEG}, formula (8.5.1):
$$
{HH}^i(\G\ltimes\sym M)=
\begin{cases}(\sym M)^\G &\text{if}\enspace i=0\\
(M^*\o\sym M)^\G&\text{if}\enspace i=1\\
\La^2M^*\o(\k[\SS]\oplus \sym M)^\G
&\text{if}\enspace i=2\\
0&\text{if}\enspace i>2.
\end{cases}
$$

It follows  by Morita equivalence
of Hochschild homology and Van den Bergh duality 
$\La^2M^*\o{HH}_\idot(\Pi)\cong  {HH}^{2-\hdot}(\Pi),$ see \cite{VdB},
that we have
\begin{eqnarray}\label{cbeg}
&\La^2M^*\o{HH}_2(\Pi)\cong {HH}^0(\Pi)\cong
 {HH}^0(\G\ltimes\sym M)\cong (\sym M)^\G,\\
&\La^2M^*\o{HH}_0(\Pi)\cong {HH}^2(\Pi)\cong {HH}^2(\G\ltimes\sym M)\cong
\La^2M^*\o(\k[\SS]\oplus \sym M)^\G.\nonumber
\end{eqnarray}

The  isomorphisms in the second line of \eqref{cbeg} imply that
$({HH}_0(\Pi))_+\cong (\sym M)^\G_+.$
Hence, in terms of the integer coefficients $\phi_k(\bc)_{oo}$
in the Hilbert series of $(\sym M)^\G$,
see Corollary \ref{b},
we obtain the formula
$$
h\big(\O(\Pi)\big)=
h\big(\sym\big({HH}_0(\Pi)_+\big)\big)=h\big(\sym\big((\sym M)^\G_+\big)\big)
=\prod_{k\geq1} \frac{1}{(1-t^k)^{\phi_k(\bc)_{oo}}}.
$$

Further, the vector space $\La^2M$ has weight degree $2$.
Hence, from
 the first line of \eqref{cbeg} we get a graded space isomorphism
${HH}_2(\Pi)\cong(\sym M)^\G\langle 2\rangle,$
where $\langle 2\rangle$ denotes grading shift by 2. Thus,
we find
$$
h\big(\sym {HH}_2(\Pi)\big)=h\big(\sym\big((\sym M)^\G\langle 2\rangle\big)\big)
=\prod_{k\geq0} \frac{1}{(1-t^{k+2})^{\phi_k(\bc)_{oo}}}.
$$

Using  the above formulas for  Hilbert series,
the last equation can be rewritten as
$$
\prod_{k\geq1} \frac{1}{(1-t^k)^{\phi_k(\bc)_{oo}}}=
\ze(Q,t)\cd\prod_{k\geq0} \frac{1}{(1-t^{k+2})^{\phi_k(\bc)_{oo}}}.
$$
Thus, we obtain
$$\ze(Q,t)=
\prod_{k\ge 1}\frac{1}{(1-t^{k})^{\phi_k(c)_{00}-\phi_{k-2}(c)_{00}}}
=\prod_{k\geq1}
\frac{1}{(1-t^{k})^{\varphi_k(\bc)_{oo}}}.
$$
This equation
is equivalent to \eqref{curious}, and we are done.
\end{proof}

\begin{remark} The identity  in \eqref{curious} may be also
verified directly, case by case,
for each Dynkin diagram of type $A,D,E.$
To do so, one has to know
explicit formulas for the determinant of
the $t$-analogue of the  Cartan matrix of the corresponding 
extended Dynkin graph.
Equivalently, for the  polynomial $D(t):=\det
(\bone-t\cdot\bc+t^2\cd\bone)$, one finds the following formulas, 
see \cite[\S3.13]{Lu}:
\begin{eqnarray*}
&A_n:     & D(t)=(1-t^{n+1})^2;\\
&D_n:     & D(t)=(1-t^4)^2(1-t^{2n-4})/(1-t^2);\\
&E_6:     & D(t)=(1-t^4)(1-t^6)^2/(1-t^2);\\
&E_7:     & D(t)=(1-t^4)(1-t^6)(1-t^8)/(1-t^2);\\
&E_8:     & D(t)=(1-t^4)(1-t^6)(1-t^{10})/(1-t^2).
\end{eqnarray*}
\end{remark}


\end{document}